\documentclass[amsmath,12pt,a4paper]{amsart}
\usepackage{amsfonts, amssymb, amsmath, amscd}
\usepackage{amsthm}
\usepackage[pdftex]{graphicx}
\usepackage{multicol}
\usepackage{algorithm}
\usepackage{algpseudocode}
\usepackage{tikz}
\usepackage{xcolor}
\usepackage{hyperref}
\usepackage{ragged2e,caption}
\usepackage{pgfplots}
\usepackage{pifont}
\pgfplotsset{compat=1.17}
\definecolor{lightblue}{rgb}{0.68,0.85,0.90} 
\usepackage{bm,bigints, comment}
\usepackage{mathrsfs}
\usepackage{color}
\usepackage{enumitem} 
\newtheorem{theorem}{Theorem}

\newtheorem{Corollary}{Corollary}

\newtheorem{lemma}{Lemma}

\newtheorem{definition}{Definition}
\newtheorem{example}{Example}
\title[The sliding tile puzzle, roots to polynomials, and $\textbf{P}$ vs. $\textbf{NP}$ complexity] {The sliding tile puzzle, roots to polynomials, and $\textbf{P}$ vs. $\textbf{NP}$ complexity}

\author{R. Burson}

\begin{document}
\maketitle
\begin{abstract}
This work explores the relationship between solution space and time complexity in the context of the $\textbf{P}$ vs. $\textbf{NP}$ problem, particularly through the lens of the sliding tile puzzle and root finding algorithms. We focus on the trade-off between finding a solution and verifying it, highlighting how understanding the structure of the solution space can inform the complexity of these problems. By examining the relationship between the number of possible configurations and the time complexity required to traverse this space we demonstrate that the minimal time to verify a solution is often smaller than the time required to discover it. Our results suggest that the efficiency of solving $\textbf{NP}$-complete problems is not only determined by the ability to find solutions but also by how effectively we can navigate and characterize the solution space. This study contributes to the ongoing discourse on computational complexity, particularly in understanding the interplay between solution space size, algorithm design, and the inherent challenges of finding versus verifying solutions.
\end{abstract}

\section{The sliding puzzle}\label{Sec1}
\subsubsection{Introduction}
The sliding tile puzzle is a popular combinatorial game invented by Noyes Chapman in 1879 and patented as the "Gem Game" or "Gem Puzzles," later widely popularized by Samuel Loyd. It consists of a grid of numbered tiles that slide within a frame, with the objective of arranging them in a specified order by sliding numbered blocks into an adjacent empty space. The solution involves a sequence of steps or moves required to reach the desired arrangement. Over the years, the puzzle has inspired both casual and competitive engagement within puzzle-solving communities, contributing to significant research in mathematical and computational areas such as graph theory, group theory, and topology \cite{BrunckFlorestan}, \cite{HasanDler}, \cite{SpitzEdward}, \cite{ArcherAaron}, \cite{SlocumJerry}. Recent studies have also examined machine learning and AI approaches to solving such puzzles, exploring their relevance to decision-making processes \cite{deOlBryan}, \cite{ChenMen}, \cite{KorfRich}, \cite{KordRich2}, \cite{Ryuhei}, \cite{Wang}. Our work focuses on the connection between the sliding tile puzzle and the famous $\textbf{P}$ vs. $\textbf{NP}$ problem \cite{GoldOded}, \cite{WigAvi}, \cite{FortnowLance} , \cite{cooks} specifically examining the computational complexity of solving particular puzzle configurations. It is well-established that some puzzle configurations, such as the $15$-puzzle has no solutions ,while others are solvable in deterministic polynomial time \cite{GozonMarc}. We aim to explore how the time complexity required to check if a solution exists correlates with the time complexity to discover the necessary moves. Similar to Gozen and Yu’s research on the complexity of solving the Generalized Shortest Traveling Path problem \cite{GozonMarc}, our investigation seeks to demonstrate that solvable sliding tile puzzles are deterministic in nature but may require varying amounts of time based on the initial configuration and especially the dimension of the puzzle (as one generalizes the problem from matrices to tensors). Any algorithm addressing a sliding puzzle must apply moves sequentially to transform the initial configuration, underscoring the deterministic nature of these puzzles. In the context of the sliding tile puzzle, minimal time complexity can be understood as the amount of time needed to apply each move sequentially. As each move reduces the distance between the current puzzle configuration and the solved state, the complexity is linked to the number of moves required to reach the solution. This insight aligns with established computational theories, particularly in the context of the $\textbf{P}$ vs. $\textbf{NP}$ problem, where the deterministic nature of solving puzzles implies that the time complexity is proportional to the number of required moves. As suggested by Gozen and Yu \cite{GozonMarc}, such complexity arises from applying operations in a systematic order, ensuring that the solution space is explored exhaustively in polynomial time.\\

\subsubsection{Preliminaries}
Unlike previous research we will examine the sliding puzzle through a matrix representation rather then using group permutations. Each puzzle configuration corresponds to a matrix with components in $\mathbb{N}\cup \{\varnothing\}$ where $\varnothing$ is a character used to denote the empty space in the matrix. Each legal move can be viewed as a transformation of the matrix by swapping an adjacent element with the empty space $\varnothing$ and there are precisely four legal transformations. \\

\begin{definition}\label{Def1}
 (Sliding tiles $\mathcal{S}$) let $n\in \mathbb{N}$ a sliding tile is an element of the set of $n\times n$ matrices with components in $\mathbb{N}\cup\{\varnothing\}$ written $M_n(\mathbb{N}\cup\{\varnothing\})$ such that no two components are equal and are less then $n$.  Each matrix element represents a numbered tile or an empty space, which allows for the tile shifting within the matrix. The set of all sliding tiles $\mathcal{S}$ is the union

\begin{equation}\label{Equation 1}
\mathcal{S} = 
\bigcup_{n=1}^\infty \textbf{T}_n
\end{equation}
 with 
 
\begin{equation}\label{Equation 2}
\textbf{T}_n = \left\{ \mathbf{A} \in M_n\left(\mathbb{N} \cup \{\varnothing\}\right) : 
\begin{array}{l}
a_{ij} \leq n-1 \text{ for all entries } a_{ij}, \\
a_{ij} \neq a_{kl} \text{ whenever } (i, j) \neq (k, l)\\
\exists ~(i,j)~~\textbf{s.t}~~a_{ij} = \varnothing
\end{array}
\right\}.
\end{equation} 

\end{definition}

\begin{lemma}\label{lem1} $|\textbf{T}_n| = n^2!$ for $n\in \mathbb{N}$. 

\begin{proof}
Consider any bijection $\mathcal{F}:~\textbf{T}_n\rightarrow \{1,2,..., n^2!\}$ then the cardinality of $\textbf{T}_n$ is provided by 

\begin{equation}\label{Equation 3}
|\textbf{T}_n| = |\mathcal{F}(\textbf{T}_n)| = |\{1,2,...,n^2!\}| = n^2!.
\end{equation}

\end{proof}
\end{lemma}

\begin{definition}\label{Def3} (Four legal transformations) let $(n,k)\in \mathbb{N}\times \{1,2,3,4\}$ and consider a sliding tile $\textbf{A}\in \mathcal{S}$. Define $\Sigma_k: \mathcal{S}\rightarrow \mathcal{S}$ using the rule 

\begin{equation}\label{Equation 4}
(\textbf{A})_{ij} \mapsto \begin{cases}
(\textbf{A})_{\rho_k(i,j)} & \text{if $(\textbf{A})_{ij} =\varnothing$}\\
(\textbf{A})_{\tau_k(i,j,\textbf{A})} & \text{if $(\textbf{A})_{ij} \neq \varnothing$}\\
\end{cases}
\end{equation}

where $\rho_k :\mathbb{N}\times
 \mathbb{N}\rightarrow \mathbb{N}\times \mathbb{N}$ is a function that provides the indices to the matrix

 \begin{equation}\label{Equation 5}
  \rho_k(i,j) = 
\begin{cases}
\biggl(i ,j+1\biggr) & \text{ if $k =1$ (up)}\\
\biggl(i , j-1\biggr) & \text{ if $k = 2$ (down)}\\
\biggl(i+1, j\biggr) & \text{ if $k = 3$ (right)}\\
\biggl(i-1 ,j \biggr) & \text{ if $k =4$ (left)}
\end{cases}
 \end{equation}
 
 and the function $\tau_k : \mathbb{N}\times \mathbb{N} \times \mathcal{S}\rightarrow \mathbb{N} \times \mathbb{N}$ is defined in terms of $\rho_k$
 
\begin{equation}\label{Equation 6}
\tau_k(i,j,\textbf{A})= \begin{cases}
 (i , j) & \textit{If $(\textbf{A})_{\lambda_k(i,j)} \neq \varnothing $} \\
\lambda_k(i,j) & \textit{If $(\textbf{A})_{\lambda_k(i,j)} = \varnothing $}  
\end{cases}
\end{equation}

with $\lambda_k: \mathbb{N}\times \mathbb{N}\rightarrow \mathbb{N}\times \mathbb{N}$ explicitly provided by the rule

\begin{equation}\label{Equation 7}
\lambda_k(i,j) = 
\begin{cases}
(i,j-1) & \text{if $\rho_k(i,j) = (i ,j+1)$}\\
( i, j+1)& \text{if $\rho_k(i,j) = (i ,j-1)$}\\
(i-1, j)& \text{if $\rho_k(i,j) = (i+1, j)$}\\
 (i+1, j )& \text{if $\rho_k(i,j) = (i-1, j)$}.\\
\end{cases}
\end{equation}

We say a transformation or legal move of a sliding tile is function $$\phi_k: \mathcal{S}\rightarrow  \mathcal{S}$$ provided as

\begin{equation}\label{Equation 8}
(\textbf{A})_{ij}\mapsto \begin{cases}
 (\textbf{A})_{ij} & \text{if $(\textbf{A})_{ij}\neq \varnothing$, $(\textbf{A})_{i(j+1)}\neq \varnothing$, $(\textbf{A})_{i(j-1)}\neq \varnothing$, $(\textbf{A})_{(i+1)j}\neq \varnothing$, $(\textbf{A})_{(i-1)j}\neq \varnothing$}\\
 \Sigma_k(\textbf{A})_{ij} & \text{Otherwise}. 
\end{cases}
\end{equation}

The unique transformations of the sliding puzzle are denoted using index notation \(\phi_1\), \(\phi_2\), \(\phi_3\), and \(\phi_4\) and they correspond to the movements up, down, right, and left, respectively. In any initial configuration these are the four unique transformations that be applied to any sliding puzzle. \\
\end{definition}

\begin{definition}\label{Def4} (Goal state) For any $n\in \mathbb{N}$ the goal state matrix is written $\textbf{A}_{\textbf{goal}}$ and defined as 

 \begin{equation}\label{Equation 9}
(\textbf{A}_{goal})_{ij} = \begin{cases}
i & \text{if $i \leq n$ and $j=1$} \\
(j-1)n+i & \text{if $i \leq n-1$ and $2\le j\le n$}\\
\varnothing & \text{if i=n=j}.
\end{cases}
\end{equation}
\end{definition}

\begin{definition}\label{Def5}(Solvable sliding puzzles) Let $n\in \mathbb{N}$ and $\textbf{A}\in M_n(\mathbb{N}\cup\{\varnothing\})$ the set of solvable sliding tiles is denoted $\mathcal{S}_n(\mathbb{N}\cup\{\varnothing\})\subset \mathcal{S}$ and is provided by  set 

\begin{equation}\label{Eq10}
\biggl\{\textbf{A}\in M_n(\mathbb{N}\cup \{\varnothing\})~:~ \exists ~\epsilon:~\mathbb{N}\rightarrow \{1,2,3,4\}~ \wedge~ k\in \mathbb{N}~~ \textbf{s.t}~~ (\phi_{\epsilon(i)}\circ \phi_{\epsilon(i)}\circ \cdots \circ \phi_{\epsilon(k)})\circ \textbf{A} = \textbf{A}_{\textbf{goal}}\biggr\}.
\end{equation}

\end{definition}

\begin{example}\label{Ex1} (Solving a sliding puzzle provided an initial configuration) Set \( n = 4 \) and consider the initial configuration  

\[
\textbf{A} = 
\begin{bmatrix}
1 & \varnothing & 2 & 4 \\
5 & 6 & 3 & 8 \\
9 & 10 & 7 & 11\\
13 & 14 & 15 & 12
\end{bmatrix}.
\]

The sequence of steps to solve the goal state $\textbf{A}_{\textbf{goal}}$ is illustrated in the table below (and is always given by a finite number of compositions of some transformation)

\[
\begin{array}{|c|c|c|l|}
\hline
\text{Step} & \text{Transformation} & \text{Resulting Matrix} & \text{Description} \\
\hline
1 & \textbf{A} \mapsto \phi_3(\textbf{A}) & 
\begin{bmatrix}
1 & 2 & \varnothing & 4 \\
5 & 6 & 3 & 8 \\
9 & 10 & 7 & 11\\
13 & 14 & 15 & 12
\end{bmatrix} & \text{Move the empty tile right.} \\
\hline
2 &  \phi_3(\textbf{A}) \mapsto \phi_2(\phi_3(\textbf{A})) & 
\begin{bmatrix}
1 & 2 & 3 & 4 \\
5 & 6 & \varnothing & 8 \\
9 & 10 & 7 & 11\\
13 & 14 & 15 & 12
\end{bmatrix} & \text{Move the empty tile down.} \\
\hline
3 &  \phi_2(\phi_3(\textbf{A})) \mapsto \phi_2(\phi_2(\phi_3(\textbf{A}))) & 
\begin{bmatrix}
1 & 2 & 3 & 4 \\
5 & 6 & 7 & 8 \\
9 & 10 & \varnothing & 11\\
13 & 14 & 15 & 12
\end{bmatrix}& \text{Move the empty tile down.} \\
\hline
4 &  \phi_2(\phi_2(\phi_3(\textbf{A}))) \mapsto \phi_3(\phi_2(\phi_2(\phi_3(\textbf{A})))) & 
\begin{bmatrix}
1 & 2 & 3 & 4 \\
5 & 6 & 7 & 8 \\
9 & 10 & 11 & \varnothing\\
13 & 14 & 15 & 12
\end{bmatrix} & \text{Move the empty tile right.} \\
\hline
5 &  \phi_3\phi_2(\phi_2(\phi_3(\textbf{A}))) \mapsto \phi_2(\phi_3(\phi_2(\phi_2(\phi_3(\textbf{A}))))) & 
\begin{bmatrix}
1 & 2 & 3 & 4 \\
5 & 6 & 7 & 8 \\
9 & 10 & 11 & 12 \\
13 & 14 & 15 & \varnothing
\end{bmatrix} & \text{Move the empty tile down.} \\
\hline
\end{array}
\]

Observe the goal state $\textbf{A}_{goal}$ was achieved using only two types of movements, namely right and down, and other transformations such as \(\phi_1\) (up) and \(\phi_4\) (left) were not needed in this specific sequence but can be applied for other scenarios. No component in the first two columns where altered they remained the same in every transformation. We recognize our steps needed to arrive at $\textbf{A}_{goal}$ can be written as finite amount of composition of functions $\{\phi_1,\phi_2, \phi_3,\phi_4\}$ and all such finite compositions. \\

\end{example}

\begin{definition}\label{Def6}  For $k \in \mathbb{N}$ define the set of possible consecutive transformations as

\begin{equation}\label{Eq11}
\mathbb{E}_k = \{(\phi_{\epsilon(1)},\phi_{\epsilon(1)},\cdots, \phi_{\epsilon(k)}) : ~~~~ \epsilon: \mathbb{N}\rightarrow \{1,2,3,4\}\}.
\end{equation}

Let $n\in \mathbb{N}$ then a general solution to the sliding puzzle game is a function $$f:~ \mathcal{S}_n(\mathbb{N}\cup\{\varnothing\})\rightarrow \bigcup_{k=1}^{\infty}{\mathbb{E}_k}$$ such that

\begin{equation}\label{Eq12}
f(\textbf{A}) \circ \textbf{A} = \textbf{A}_{goal}, \indent \forall \textbf{A}\in \mathcal{S}_n(\mathbb{N}\cup\{\varnothing\}).
\end{equation}

\end{definition}

\begin{definition}\label{Def7}
Let $n\in \mathbb{N}$ and consider a sliding tile $\textbf{A}\in \mathcal{S}_n(\mathbb{N}\cup\{\varnothing\}).$  We say the minimal number of \textit{moves} to solve the sliding tile is defined as the mapping $\Psi:~\mathcal{S}_n(\mathbb{N}\cup \{\varnothing\})\rightarrow \mathbb{N}$ provided by the rule $$\textbf{A} \mapsto \textbf{Min}\biggl\{k:~ \biggl(\phi_{\epsilon(1)}\circ \phi_{\epsilon(2)}\circ \cdots \circ \phi_{\epsilon(k)}\biggr)\circ \textbf{A} = \textbf{A}_{goal},~~ \forall~~ \epsilon:~\mathbb{N}\rightarrow \{1,2,3,4\}\biggr\}.$$
\end{definition}

\subsubsection{Elementary Results}
In this section we explore some basic fundamental results concerning the sliding puzzle game. We start with the following lemma

\begin{lemma}\label{Lem2}
For $\textbf{A}\in \mathcal{S}_n(\mathbb{N}\cup \{\varnothing\})$ we have 

\begin{equation}\label{Eq13}
\phi_1(\phi_2 (\textbf{A}))= \textbf{A} ,\indent \phi_2(\phi_1 (\textbf{A}))= \textbf{A} 
\end{equation}

and 

\begin{equation}\label{Eq14}
\phi_4(\phi_3 (\textbf{A}))= \textbf{A} ,\indent \phi_4(\phi_3 (\textbf{A}))= \textbf{A} 
\end{equation}

\end{lemma}

\begin{theorem}\label{Thm1}(Goal state Equivalence)

\begin{equation}\label{Eq15}
\textbf{A}\in \mathcal{S}_n(\mathbb{N}\cup \{\varnothing\})\Longleftrightarrow \exists k\in \mathbb{N} \wedge  \epsilon:~\mathbb{N}\rightarrow \{1,2,3,4\}, ~~ \biggl(\phi_{\epsilon(1)}\circ \phi_{\epsilon(2)} \circ \cdots \phi_{\epsilon(k)}\biggr)\circ \textbf{A}_{goal} = \textbf{A}
\end{equation}

\begin{proof}
Let $n\in\mathbb{N}$ and assume $\textbf{A}\in \mathcal{S}_n(\mathbb{N}\cup \{\varnothing\})$ by definition there exist a value $k\in \mathbb{N}$ and a function $f:~\mathbb{N}\rightarrow \{1,2,3,4\}$ such that $$\biggl(\phi_{f(1)}\circ \phi_{f(2)} \circ \cdots \phi_{f(k)}\biggr)\circ \textbf{A}  =  \textbf{A}_{goal}.$$ Define the function $\epsilon:~\mathbb{N}\rightarrow \{1,2,3,4\}$ using the rule

 $$j\mapsto \begin{cases}
1 & \text{if $f(k-j+1) = 2$}.\\
2 & \text{if $f(k-j+1) = 1$}.\\
3 & \text{if $f(k-j+1) = 4$}.\\
4 & \text{if $f(k-j+1) = 3$}.
 \end{cases}$$

  Then 

\begin{equation}\label{Eq16}
\begin{split}
\biggl(\phi_{f(1)}\circ \phi_{f(2)} \circ \cdots \phi_{f(k)}\circ \textbf{A} \biggr) &=  \textbf{A}_{goal}.\\
& \Updownarrow\\
\biggl(\phi_{\epsilon(1)}\circ \phi_{\epsilon(2)} \circ \cdots \phi_{\epsilon(k)}\biggr) \circ \biggl(\phi_{f(1)}\circ \phi_{f(2)} \circ \cdots \phi_{f(k)}\circ \textbf{A} \biggr) &=  \biggl(\phi_{\epsilon(1)}\circ \phi_{\epsilon(2)} \circ \cdots \phi_{\epsilon(k)}\biggr)\circ  \textbf{A}_{goal}.\\
& \Updownarrow\\
\textbf{A}  &= \biggl(\phi_{\epsilon(1)}\circ \phi_{\epsilon(2)} \circ \cdots \phi_{\epsilon(k)}\biggr)\circ \textbf{A}_{goal} 
\end{split}
\end{equation}
\end{proof}
\end{theorem}

Theorem \ref{Thm1} shows that instead of looking for solutions to the sliding puzzle game with an initial configuration one can instead look for a solution starting at the goal state that reaches the initial configuration. As in Example \ref{Ex1} we have shown $\phi_2(\phi_3(\phi_2(\phi_2(\phi_3(\textbf{A}))))) = \textbf{A}_{goal}$ and hence from Theorem \ref{Thm1} we immediately must have $\phi_3(\phi_2(\phi_2(\phi_3(\phi_2(\textbf{A}_{goal}))))) = \textbf{A}$.

\begin{theorem}\label{Thm2}
Let $n\in\mathbb{N}$ and take $\textbf{A}\in \mathbb{S}_n(\mathbb{N}\cup \{\varnothing\})$. Then 

\begin{equation} \label{Eq17}
\Psi(\textbf{A})\le \log_4(n^2!).
\end{equation}

\begin{proof}
The proof is rather instructive. First assume $\textbf{A}\in \mathbb{S}_n(\mathbb{N}\cup \{\varnothing\})$ and then by construction there is a function $\epsilon:~\mathbb{N}\rightarrow \{1,2,3,4\}$ such that $$\biggl(\phi_{\epsilon(1)}\circ \phi_{\epsilon(2)} \circ \cdots \phi_{\epsilon(k)}\circ \textbf{A} \biggr) =  \textbf{A}_{goal}, \indent k = \Psi(\textbf{A}).$$

Each transformation $\{\phi_{\epsilon(1)},\phi_{\epsilon(2)}, \cdots, \phi_{\epsilon(k)} \}$ produces a new element in the set $\mathcal{S}_n(\mathbb{N}\cup \{\varnothing\})$ and because $k = \Psi(\textbf{A})$ is the minimal value that produces the goal state $\textbf{A}_{goal}$. For each transformation there is a maximum of $4$ different possibilities they can be particularly $\{\phi_1,\phi_2,\phi_3,\phi_4\}$. However, regardless of the number of transformations we consecutively apply we cannot exceed the number of elements in the set $M_n(\mathbb{N}\cup\{\varnothing\})$ as $\mathcal{S}_n(\mathbb{N}\cup\{\varnothing\})\subset M_n(\mathbb{N}\cup\{\varnothing\})$.  Therefore, we must have 

\begin{equation}\label{Eq18}
4^{\Psi(\textbf{A})} \le |M_n(\mathbb{N}\cup\{\varnothing\})| = n^2!
\end{equation}

which gives us

\begin{equation}\label{Eq19}
\Psi(\textbf{A})\le \log_4(n^2!).
\end{equation}
\end{proof}
\end{theorem}

\begin{theorem}\label{Thm3}
Let $n\in \mathbb{N}$ with $n>1$ then $$|\mathcal{S}_n(\mathbb{N}\cup \{\varnothing\})| \le 4 \biggl(3^{\lfloor \dfrac{\log_4(n^2!)-1}{2}\rfloor }\biggr)\biggl(4^{\lfloor \dfrac{ \log_4(n^2!)-1}{2}}\rfloor\biggr)+4.$$ 

\begin{proof}
 The general solution $f(\textbf{A})$ is a sequence of moves $\phi_{\epsilon(1)}\circ \phi_{\epsilon(2)}\circ \cdots \circ \phi_{\epsilon(k)}$ with $\epsilon:\{1,2,\cdots \Psi(\textbf{A})\}\rightarrow \{1,2,3,4\}$ and this function exist because $\Psi(\textbf{A})\ge 1$ meaning at one least one move must be applied to reach the goal state matrix. The maximum different possible choices of $\epsilon(\Psi(\textbf{A}))$ is $four$ because $|\epsilon(\{1,2,\cdots \Psi(\textbf{A})\})| = |\{1,2,3,4\}| = 4$,  but then $\epsilon_{(\Psi(\textbf{A})-1)} $  cannot take on the value $i$ that undoes the operator $\epsilon (\psi(\textbf{A}))$ as we see from Lemma \ref{Lem2}. Consider the set $$\{\textbf{A}\in \mathcal{S}_{n}(\mathbb{N}\cup\{\varnothing\}):~\Psi(\textbf{A})\ge 2\}.$$ These are the puzzles solvable in exactly two moves or more. Allow $\epsilon: \mathbb{N} \rightarrow \mathbb{N}$ to be any sequence that solves a sliding tile in this particular set. Then exactly half the choices of each $\epsilon(k)$ have three possibilities and the other half has four choices.  We also notice that $$|\{\textbf{A}\in \mathcal{S}_{n}(\mathbb{N}\cup\{\varnothing\}):~\Psi(\textbf{A})=1\}|\le |\{\phi_1,\phi_2,\phi_3,\phi_4\}| $$ because any tile that can be solved in one move can be achieved in one move from the goal equiv lance theorem \ref{Thm1}. If it takes one move to go from $\textbf{A}$ to the matrix $\textbf{A}_{goal}$ then it takes one move to go from $\textbf{A}_{goal}$ to $\textbf{A}$ but there is also four possible matrices that this can be $$\biggl\{\phi_1(\textbf{A}_{goal}),  \phi_2(\textbf{A}_{goal}), \phi_3(\textbf{A}_{goal}), \phi_4(\textbf{A}_{goal})\biggr\}.$$ (we will ignore the fact that we cannot swap the empty tile down or to the right starting at the goal matrix). Set

 $$M = \textbf{Max}\left(\{ \Psi(\textbf{A})-1:~~~\textbf{A}\in \mathcal{S}_n(\mathbb{N}\cup\{\varnothing\}\})\right)$$  
 
We have $$|\{\textbf{A}\in \mathcal{S}_{n}(\mathbb{N}\cup\{\varnothing\}):~\Psi(\textbf{A})\ge 2\}|\le 4\biggl(3^{\lfloor\dfrac{M}{2}\rfloor}\biggr)\biggl(4^{\lfloor\dfrac{M}{2}\rfloor}\biggr).$$ Here $\lfloor \cdot \rfloor$ denotes the floor function. From Theorem \ref{Thm2} we know $M\le \log_4(n^2!)-1$ then we conclude with the following inequalities

\begin{equation*}
\begin{split}
|\mathcal{S}_n(\mathbb{N}\cup \{\varnothing\})| & \le |\{\textbf{A}\in \mathcal{S}_{n}(\mathbb{N}\cup\{\varnothing\}):~\Psi(\textbf{A})\ge 2\}|+|\{\textbf{A}\in \mathcal{S}_{n}(\mathbb{N}\cup\{\varnothing\}):~\Psi(\textbf{A})=1\}|\\
&  \le 4 \biggl(3^{\lfloor\dfrac{M}{2}\rfloor}\biggr)\biggl(4^{\lfloor\dfrac{M}{2}\rfloor}\biggr)+|\{\textbf{A}\in \mathcal{S}_{n}(\mathbb{N}\cup\{\varnothing\}):~\Psi(\textbf{A})=1\}|\\
& \le 4 \biggl(3^{\lfloor \dfrac{\log_4(n^2!)-1}{2}\rfloor }\biggr)\biggl(4^{\lfloor \dfrac{ \log_4(n^2!)-1}{2}}\rfloor\biggr)+ |\{\phi_1,\phi_2,\phi_3,\phi_4\}|\\
& \le  4 \biggl(3^{\lfloor \dfrac{\log_4(n^2!)-1}{2}\rfloor }\biggr)\biggl(4^{\lfloor \dfrac{ \log_4(n^2!)-1}{2}}\rfloor\biggr)+4
 \end{split}
\end{equation*}

\end{proof}
\end{theorem}

\begin{theorem}\label{Thm4}
If $n\ge 3$ then 

$$|\mathcal{S}_n(\mathbb{N}\cup \{\varnothing\})| \le 4(n^2-n-4) $$

\begin{proof}
The proof is combinatorial we consider the character $\varnothing$ moving along the grid in ever possible component and count how many matrix we can produce from any application $\{\phi_1,\phi_2,\phi_3,\phi_4\}$ see figure 1 below. let $\textbf{A}\in \mathcal{S}_n(\mathbb{N}\cup\{\varnothing\})$ with $n>3$. We notice for each possible $(i,j)$ such that $a_{ij} = \varnothing$ there is either $2$, $3$, or $4$ different possible application we can apply to matrix configuration to move the tile $\varnothing$. If $a_{11} = \varnothing$, $a_{1n} = \varnothing$, $a_{n1}= \varnothing$, or $a_{nn} = \varnothing$ then we can only apply two possible transformations in such case (these are the four corners of the matrix). We cannot apply all of the four transformations but each time we apply any one them we get a new matrix. There are four edges to the matrix then there are $4*2 = 8$ possible different matrix from such application. Next for any column or row in the matrix there is exactly $n$ different elements inside them. If we neglect the four corners there are $(n-2)$ elements that have precisely three possible transformation we can apply to them (i.e the first row minus the two corners in the first row). Consider the specific components along the edges in the following sets

 $$\mathbb{W}_1 = \{a_{ij}:~ i=1, 1<j< n\},$$ 

$$\mathbb{W}_2 = \{a_{ij}:~ j=n, 1<i< n\},$$

$$\mathbb{W}_3 = \{a_{ij}:~ j=1, 1<i< n\},$$

$$\mathbb{W}_4 = \{a_{ij}:~ i=n, 1<j< n\}.$$

If $\varnothing \in \mathbb{W}_1\cup \mathbb{W}_2\cup\mathbb{W}_3 \cup \mathbb{W}_4$ then we can only apply three transformation to any such matrix to get a new resulting matrix. Each set contains $(n-2)$ elements and so there are $3(n-2)$ maximum possible matrix we can produce from applying any of the three transformation to the $(n-2)$ different components, but again there are four sets so there is $4*3(n-2) = 12(n-2)$ different matrix we can move transform three different ways. Now for any component $a_{ij}$ not along the four outside edges there is a maximum of four different application we can apply to attempt to reach the goal state. To get this value we subtract the number of components along the edges from the total number of components to the matrix which is $n^2- 4n$. So there are no more then $4(n^2- 4n)$ matrices we can get from applying any transformation in this case. Therefore, 

$$|\mathcal{S}_n(\mathbb{N}\cup \{\varnothing\})| \le 8 + 12(n-2) + 4(n^2-4n) = 4(n^2-n-4)$$

\begin{center}
\begin{tikzpicture}
    \def\n{5} 
    
    \foreach \i in {1,...,\n} {
        \foreach \j in {1,...,\n} {
            \pgfmathsetmacro\x{(\j-1) + 0.5}
            \pgfmathsetmacro\y{(-\i+1) - 0.5}
            
            \ifnum\i=1
                \ifnum\j=1
                    \fill[red] (\j-1,-\i+1) rectangle (\j,-\i); 
                    \node at (\x,\y) {2}; 
                \else
                    \ifnum\j=\n
                        \fill[red] (\j-1,-\i+1) rectangle (\j,-\i); 
                        \node at (\x,\y) {2}; 
                    \else
                        \fill[green] (\j-1,-\i+1) rectangle (\j,-\i); 
                        \node at (\x,\y) {3}; 
                    \fi
                \fi
            \else
                \ifnum\i=\n
                    \ifnum\j=1
                        \fill[red] (\j-1,-\i+1) rectangle (\j,-\i); 
                        \node at (\x,\y) {2}; 
                    \else
                        \ifnum\j=\n
                            \fill[red] (\j-1,-\i+1) rectangle (\j,-\i); 
                            \node at (\x,\y) {2}; 
                        \else
                            \fill[green] (\j-1,-\i+1) rectangle (\j,-\i); 
                            \node at (\x,\y) {3}; 
                        \fi
                    \fi
                \else
                    \ifnum\j=1
                        \fill[green] (\j-1,-\i+1) rectangle (\j,-\i); 
                        \node at (\x,\y) {3}; 
                    \else
                        \ifnum\j=\n
                            \fill[green] (\j-1,-\i+1) rectangle (\j,-\i); 
                            \node at (\x,\y) {3}; 
                        \else
                            \fill[lightblue] (\j-1,-\i+1) rectangle (\j,-\i);
                            \node at (\x,\y) {4}; 
                        \fi
                    \fi
                \fi
            \fi
        }
    }

    \foreach \i in {0,...,\n} {
        \draw (\i,0) -- (\i,-\n); 
        \draw (0,-\i) -- (\n,-\i); 
    }

    \node at (-0.5,0.5) {\(a_{11}\)};
    \node at (\n+0.5,0.5) {\(a_{1n}\)};
    \node at (-0.5,-\n+0.5) {\(a_{n1}\)};
    \node at (\n+0.5,-\n+0.5) {\(a_{nn}\)};
\end{tikzpicture}
\end{center}
    \vspace{0.5cm}
    \textbf{Figure 1:} An example of the number of moves in a $5x5$ sliding tile along every component. 
    \label{fig1}
    
\end{proof}    
\end{theorem}

\begin{Corollary}\label{Cor1}
Let $n\ge 3$ for $\textbf{A}\in \mathcal{S}_n(\mathbb{N}\cup\{\varnothing \})$ $$\Psi(\textbf{A})\le 4(n^2-n-2)$$

\begin{proof}
We know $\Psi(\textbf{A})\le |\mathcal{S}_n(\mathbb{N}\cup\{\varnothing \})|$ as it cannot take more applications then the cardinality of the set itself to reach the goal state. Next, from Theorem \ref{Thm5} we know $|\mathcal{S}_n(\mathbb{N}\cup\{\varnothing \})|\le 4(n^2-n-2)$ and the result follows from the transitive property. 
\end{proof}
\end{Corollary}

\subsubsection{Time Complexity Independence}
In this short section we will discuss some topics involving the time complexity of the sliding puzzle through the perspective of the $\textbf{P}$ vs. $\textbf{NP}$ problem and provide more formal details in Section \ref{sec3}. The time complexity for the sliding puzzle is known to be at most exponential in the worst case. However, due to the existence of a "general solution" (see definition \ref{Def6}) the time to verify a solution is much less then to find it by directly writing a program with this map. Although, one does not have a explicit formula for a general solution the existence will suffice, though is not necessarily the worst possible algorithm since there is a possibility that we can write a faster algorithm in particularly one that doesn't exhaust the solution space. The time it takes to solve a given initial configuration is at most the amount of time it takes to apply the general solution (because one does not need to check the general solution and we must apply it to get to the end goal), but the time it takes to verify or check that a given solution is correct is equivalent to how long it takes to apply the arbitrary transformations but additional check that it produces the goal state $\textbf{A}_{goal}$.  Regardless of the intelligence level of a given programmer the general solution exist. The problem can be stated as followed: \\

\textbf{Problem 1:} For $\textbf{A}\in\mathcal{S}_n(\mathbb{N}\cup\{\varnothing\})$ what is the minimal amount of time that is required to find the sequence of transformations that transform $\textbf{A}$ into its goal state $\textbf{A}_{goal}$? And secondly, given a solution sequence , say $\{\phi_{\epsilon(1)},\phi_{\epsilon(2)}, \cdots, \phi_{\epsilon(k)}\}$, what is the time it takes to check that

\begin{equation}\label{Eq20}
\biggl(\phi_{\epsilon(1)}\circ \phi_{\epsilon(2)}\circ \cdots \circ \phi_{\epsilon(k)}\biggr)\circ \textbf{A} = \textbf{A}_{goal}
\end{equation}

is a true statement? \\

The time it takes to check the validity of Eq.(\ref{Eq20}) does not depend on the time it takes to find the sequence $\{\phi_{\epsilon(1)},\phi_{\epsilon(2)}, \cdots, \phi_{\epsilon(k)}\}$ but rather how long it takes to apply this sequence to a given configuration (which is nearly instantaneous provided its definition) and then after we apply them we must check that the resulting matrix $\biggl(\phi_{\epsilon(1)}\circ\phi_{\epsilon(2)}\circ \cdots \phi_{\epsilon(k)}\biggr) \circ \textbf{A} $ is indeed the goal state $\textbf{A}_{goal}$. In general to check that any two matrix are equal we must check every component. The difficulty in finding a general solution to the sliding puzzle problem arises from finding the particular value $\Psi(\textbf{A})$ (the minimum number of consecutive applications of the transformations $\{\phi_1,\phi_2,\phi_3,\phi_4\}$ to reach the goal state $\textbf{A}_{goal}$) and a corresponding sequence $\epsilon: \mathbb{N}\rightarrow \{1,2,3,4\}$ which produces a solution. If one can find a general solution $f$ to the sliding puzzle game (see definition \ref{Def6}) then the time it takes to find a solution for any given sliding tile is at at most the time it takes to apply the general solution $f$ since it might be possible to reduce the calculation of $f$ through various other methods. Reductions method may exist regardless of the representation of the map. The general solution $f$ is always correct in the sense that it always maps one to the goal state $\textbf{A}_{goal}$ one does not need to individually check the components of the matrix $f(\textbf{A})$. Nonetheless, it does take time to compute $f(\textbf{A})$ regardless of its complexity or how we may be able to represent it. One must ask why this is and the reasoning is strait forward. To check if $\biggl(\phi_{\epsilon(1)}\circ\phi_{\epsilon(2)}\circ \cdots \phi_{\epsilon(k)}\biggr) \circ \textbf{A} = \textbf{A}_{goal}$ one can apply each transformation to the matrix $\textbf{A}$ one at a time and after check if the result is equal to $\textbf{A}_{goal}$ which can be done component wise by checking the validity of the statement 

\begin{equation}\label{Eq21}
\biggl(\biggl(\phi_{\epsilon(1)}\circ\phi_{\epsilon(2)}\circ \cdots \phi_{\epsilon(k)}\biggr) \circ \textbf{A}\biggr)_{ij}  = \begin{cases}
i & \text{if $i \leq n$ and $j=1$} \\
(j-1)n+i & \text{if $i \leq n-1$ and $2\le j\le n$}\\
\varnothing & \text{if i=n=j}
\end{cases}
\end{equation}

for each $1\le i\le n$ and $1\le j\le n$. However, this is actually not the only method to validate if a given sequence of transformations produces the goal state. Instead of checking the components of the resulting matrix one can also compare the sequence $\biggl\{\phi_{\epsilon(1)},\phi_{\epsilon(2)}, \cdots \phi_{\epsilon(k)}\biggr\}$ directly to the general solution $f(\textbf{A})$ for any $\textbf{A}\in \mathcal{S}_n(\mathbb{N}\cup\{\varnothing\})$, but this requires one to know the physical representation of $f$ at $\textbf{A}$ and this generally speaking should take more time then checking $n^2$ different components of the matrix $\biggl(\phi_{\epsilon(1)}\circ\phi_{\epsilon(2)}\circ \cdots \phi_{\epsilon(k)}\biggr) \circ \textbf{A}.$ Thus using Eq.(\ref{Eq21}) is preferred for the checking process. In higher dimensions (i.e for tensors of order $3$ or more) the time to check that a given solution is correct will generally take longer the larger the dimension is because there are more components to verify and more direction to move not just left, right, up , or down. In three dimensions, there are six moves not four. While the dimension of the puzzle increases the solution space becomes exhaustively more and more complex containing multiple more elements then lower dimensional puzzles. The number of moves to any sliding tile is twice the dimension or base of the set $M_n(\mathbb{R})$ the matrices with coefficients in $\mathbb{R}$. However, regardless of the dimension of a given puzzle configuration it appears to take relatively the same time to check the solution of puzzle regardless of the dimension and so in theory this reveals that the time to find a solution to a sliding tile puzzle problem generally gets larger the higher the dimension gets. The processes of finding and verifying are independent in most cases because the way we check a problem is correct is often much different then how apply techniques to arrive at a solution. The issue relative to the $\textbf{P}$ vs. $\textbf{NP}$ problem is the time it takes computer code to get the result $f(\textbf{A})$ generally will be slower then checking if $f(\textbf{A})$ is actually a solution provided it before hand. This brings us to an important obstacle which is our ability to define a mapping that lets us compute $f$ in the minimal time, which we will need to define properly. No matter what $f$ is there may exist a $f^\prime$ where $f^\prime(\textbf{A})=f(\textbf{A}) = \textbf{A}_{goal}$  but $\mathcal{T}_M(f^\prime(\textbf{A}))\le\mathcal{T}_M(f(\textbf{A}))$  where we define $\mathcal{T}_M$ as the time it takes any Turing machine $M$ to compute a function $f$. We will give a formal definition for this function $\mathcal{T}_\mathcal{M}$ in Section \ref{sec3} when we discuss in depth the $\textbf{P}$ vs. $\textbf{NP}$ problem and formally define it but first we will discuss the problem concerning roots to polynomials to further grasp the solution space we are dealing with.

\section{Roots of Polynomials}\label{Sec2}
\subsubsection{Introduction} Musa al-Khwarizmi, a Persian mathematician, was one of the first known individuals to systematically address the roots of polynomials in his work \textit{Al-Kitab al-Mukhtasar fi Hisab al-Jabr wal-Muqabala} ("The Compendious Book on Calculation by Completion and Balancing") \cite{Rosen}, written around 780–850. This foundational text profoundly influenced both Islamic and European mathematics. In the 16th century, mathematicians such as \textit{Niccolò Fontana Tartaglia} and \textit{Gerolamo Cardano} advanced the understanding of polynomial equations, particularly cubic and quartic cases. Cardano’s publication \cite{Cardano} was pivotal in formalizing these solutions and marked a significant milestone in algebra.  By the turn of the 18th century the first known rigorous proof of the fundamental theorem of algebra was provided by Gauss \cite{Gauss}. Galois \cite{Galois2},\cite{Galois} developed what is now known as Galois Theory, which investigates the relationships between the roots of a polynomial and the symmetries (or permutations) of those roots. His work demonstrated the conditions under which polynomial equations of degree five or higher cannot be solved using radicals. Galois' key contributions laid the foundation for modern group theory and field theory, essential tools in abstract algebra. Abel \cite{Abel} also tackled the problem of solving polynomial equations of degree five or higher proving that there is no general solution in radicals for these equations, a result now known as the \textit{Abel–Ruffini Theorem} \cite{Pierpont}. Abel's work was a critical step in the development of modern algebra and served as a precursor to Galois' more comprehensive theory. Currently, a general formula that exhibits the roots to an arbitrary polynomial is not known to exist. Instead of focusing on finding general formulas our work is concerned with the time complexity of finding the roots versus checking them. Starting we will show there exist a finite set of algebraic systems to solve that exhibits the roots to any polynomial over the real numbers. First furnishing the main theorem and three examples before exploring the time complexity of the algebraic systems we present.  Further, we introduce a root finding algorithm by exhausting all a polynomials factorization equipped with Theorem \ref{Thm5} and illustrated in example \ref{Ex4}. We classify the roots of a polynomial by studying all nicely factored polynomials with less degree from which we are able to effectively determine if a polynomial has roots in a given ring $\mathcal{R}$. To avoid steering to far off topic we will gather the complexity results into the last section when we introduce the $\textbf{P}$ vs $\textbf{NP}$ problem and how it related to the sliding puzzles and the roots of polynomials. \\


\subsubsection{Preliminaries}
Recall from modern algebra if $\mathcal{R}$ is a integral domain then the ring of polynomials over $\mathcal{R}$ is defined as the set 

$$\mathcal{R}[x]=\biggl\{x\mapsto \sum_{i=0}^{d}{f(i)x^i}, ~~f:\mathbb{N}\cup\{0\}\rightarrow \mathcal{R},~ d\in \mathbb{N}\cup\{0\}\biggr\}.$$

Here $f$ is the function that represents the polynomials coefficients and $d$ denotes its degree.  \\

\begin{definition}\label{Def7}
The general solution to the roots of polynomials over an integral domain $\mathcal{R}$ is a function $\mathcal{G}:\mathcal{R}[x]\rightarrow \mathcal{R}$ defined by the rule $$ p\mapsto \{a\in \mathcal{R}:~ p(a) = 0\}.$$ For any $p\in \mathcal{R}[x]$ the general solution satisfies $p(\mathcal{G}(p)) = \{0\}$. That is if we evaluate $p$ along points in $\mathcal{G}(p)$ we get the singleton set with the zero element in the ring $\mathcal{R}$.
\end{definition}

\begin{definition}\label{Def8}
Let $\mathcal{R}$ be a integral domain. The function $\tau: \mathcal{R}[x]\rightarrow \mathbb{N}$ defined the by the rule 

$$p \mapsto |\{a\in \mathcal{R}: ~~p(a) = 0\}|$$
is called the root counting function over $\mathcal{R}$. 
\end{definition}

\begin{definition}\label{Def9}
Let $\mathcal{R}$ be a integral domain the set of irreducible polynomials over $\mathcal{R}$ is the set $\mathcal{I}_{\mathcal{R}}[x]$ defined as

$$\mathcal{I}_{\mathcal{R}}[x]= \biggl\{p\in \mathcal{R}[x]: ~ \textbf{Deg}(p)\ge 1,~~q\not\vert p,~~ \forall q \in \mathcal{R}[x]\biggr\}.$$
\end{definition}

\begin{definition}\label{Def11}
If \( \mathcal{R} \) is an integral domain, then the \textbf{absolute value} \( |\cdot| \) is a function \( |\cdot|: \mathcal{R} \rightarrow \mathbb{R} \) satisfying the following properties:
\begin{enumerate}[label=(\roman*)]
    \item \textbf{Non-negativity:} \( |x| \geq 0 \) for all \( x \in \mathcal{R} \).
    \item \textbf{Positive definiteness:} \( |x| = 0 \iff x = 0 \).
    \item \textbf{Multiplicativity:} \( |xy| = |x| \cdot |y| \) for all \( x, y \in \mathcal{R} \).
    \item \textbf{Triangle inequality:} \( |x + y| \leq |x| + |y| \) for all \( x, y \in \mathcal{R} \).
\end{enumerate}

\end{definition}

\begin{theorem}\label{Thm5}
If there exist an absolute value $|\cdot|~:\mathcal{R}\rightarrow \mathbb{R}$ then there then exist a absolute value $|\cdot|:~\mathcal{R}[x]\rightarrow \mathbb{R}$ .

\begin{proof}
Take any polynomial \( f(x) = a_0 + a_1x + \dots + a_nx^n \in \mathcal{R}[x] \), and define the function \( |\cdot| \) by the rule

\[
|f(x)| = \max\{|a_0|, |a_1|, \dots, |a_n|\},
\]
where \( |\cdot| \) is the absolute value on \( \mathcal{R} \) 

We verify that \( |\cdot| \) satisfies the axioms of an absolute value:

\begin{enumerate}[label=(\roman*)]
    \item \textbf{Non-negativity:} 
    By definition, \( |f(x)| = \max\{|a_0|, |a_1|, \dots, |a_n|\} \), and since \( |\cdot| \) on \( \mathcal{R} \) is non-negative, it follows that \( |f(x)| \geq 0 \).

    \item \textbf{Positive definiteness:} 
    If \( |f(x)| = 0 \), then \( \max\{|a_0|, |a_1|, \dots, |a_n|\} = 0 \). Since \( |\cdot| \) on \( \mathcal{R} \) satisfies positive definiteness, each \( a_i = 0 \). Thus, \( f(x) = 0 \). Conversely, if \( f(x) = 0 \), then all \( a_i = 0 \), so \( |f(x)| = 0 \).

    \item \textbf{Multiplicativity:} 
    Let \( f(x) = \sum_{i=0}^n a_i x^i \) and \( g(x) = \sum_{j=0}^m b_j x^j \). Then,
    \[
    f(x)g(x) = \sum_{k=0}^{n+m} c_k x^k,
    \]
    where \( c_k = \sum_{i+j=k} a_i b_j \). By the definition of \( |\cdot| \),
    \[
    |f(x)g(x)| = \max\{|c_0|, |c_1|, \dots, |c_{n+m}|\}.
    \]
    Using the triangle inequality and multiplicativity of \( |\cdot| \) on \( \mathcal{R} \),
    \[
    |c_k| \leq \max\{|a_i| |b_j| : i+j=k\} \leq \max\{|a_i|\} \cdot \max\{|b_j|\}.
    \]
    Thus, \( |f(x)g(x)| \leq |f(x)| \cdot |g(x)| \). Equality holds for terms achieving the maximum, so \( |f(x)g(x)| = |f(x)| \cdot |g(x)| \).

    \item \textbf{Triangle inequality:} 
    Let \( f(x), g(x) \in \mathcal{R}[x] \). Then,
    \[
    |f(x) + g(x)| = \max\{|a_i + b_i| : 0 \leq i \leq \max(n, m)\},
    \]
    where \( a_i, b_i \) are coefficients of \( f(x) \) and \( g(x) \), respectively. Using the triangle inequality on \( |\cdot| \) for \( \mathcal{R} \),
    \[
    |a_i + b_i| \leq |a_i| + |b_i|.
    \]
    Thus,
    \[
    |f(x) + g(x)| \leq \max\{|a_i| + |b_i|\} \leq \max\{|a_i|\} + \max\{|b_i|\} = |f(x)| + |g(x)|.
    \]
\end{enumerate}

Since \( |\cdot| \) satisfies all the axioms, it is an absolute value on \( \mathcal{R}[x] \).
\end{proof}

\end{theorem}

\begin{definition}\label{Def12}
Let $\mathcal{R}$ be an integral domain  in which there exist an absolute value function $|\cdot |:\mathcal{R}\rightarrow \mathbb{R}$ and consider the set of polynomials over $\mathcal{R}[x]$. We define the function $\mathcal{V}:\mathcal{R}[x] \setminus \mathcal{I}_{\mathcal{R}}[x] \to \mathcal{R}[x] \setminus \mathcal{I}_{\mathcal{R}}[x]$ by the rule

$$
p(x) \mapsto \left\{ g \in \mathcal{R}[x] \setminus \mathcal{I}_{\mathcal{R}}[x] : 
\begin{aligned}
    & g \mid p, \\
    & \deg(g) \ge \deg(q),~ \wedge~~ |g| \ge |q|,\quad \forall q \mid p
\end{aligned}
\right\}.
$$

 $\mathcal{V}(p(x))$ produces a singleton set that contains the reducible divisor of $p(x)$ which has the maximal degree among all reducible divisors of $p(x)$, but also has the largest absolute value.
\end{definition}

\begin{definition}(Root multiplicity)\label{Def13}
Let $\mathcal{R}$ be an integral domain the function $$\Upsilon : (\mathcal{R}[x]\setminus \mathcal{I}_{\mathcal{R}}[x])\times \bigcup_{p\in \mathcal{R}[x]}^{}{\{x\in \mathcal{R}:~ p(x) = 0\}}\rightarrow \mathbb{N}$$ defined by the rule 

$$(p,r)\mapsto \textbf{Max}\{m\in \mathbb{N}: ~(x-r)^m~\vert ~p(x)\}$$ is called the multiplicity of $r$ relative to $p$.  It determines how many times a root occurs in a polynomials factorization.
\end{definition}

\begin{definition}(Nicely factored polynomials)\label{Def14}
let $\mathcal{R}$ be an integral domain with an absolute value function. The set of \textit{Nicely} factored polynomials over $\mathcal{R}$ is the set $$\mathcal{N}_{\mathcal{R}}[x]= \biggl\{p\in \mathcal{R}[x]:~ p(x) \in \mathcal{V}(p(x))\biggr\}.$$ If an absolute value function $|\cdot |$ does not exist on $\mathcal{R}$ then instead we alternatively define $\mathcal{N}_{\mathcal{R}}[x]$ as the set of all reducible polynomials that are not divisible by any irreducible polynomial. So that its factorization does not depend on any irreducible. 
\end{definition}

\begin{lemma}\label{Lem3}
let $\mathcal{R}$ be an integral domain equipped with an absolute value function $|\cdot|$ and take $p\in \mathcal{N}_{\mathcal{R}}[x]$, set $k=\tau(p)$, and consider any bijection $r: \{1,2,\cdots, k\}\rightarrow\{x\in \mathcal{R}:~ p(x) = 0\}$ then there exist a function $m :\{1,2,\cdots, k\}\rightarrow \mathbb{N}$, $r :\{1,2,\cdots, k\}\rightarrow \mathbb{N}$ , and a value $c\in \mathcal{R}$ such that $$p(x) = c\prod_{i=1}^{k}{\biggl(x-r(i)\biggr)^{m(i)}},\indent \sum_{i=1}^{k}{m(i)}=\textbf{Deg}(p).$$ We will often use indicial notation and write $r_i\equiv r(i)$ and $m_i\equiv m(i)$ for any sequences.

\begin{proof}
Let $p\in \mathcal{N}_{\mathcal{R}}[x]$. First notice $\mathcal{N}_{\mathcal{R}}[x]\subset \mathcal{R}[x]\setminus \mathcal{I}_{\mathcal{R}}[x]$ and thus $p$ is reducible by construction. Since, $p$ is reducible there exist a function $m :\{1,2,\cdots, k\}\rightarrow \mathbb{N}$ and a value $c\in \mathcal{R}$ such that $$p(x) = c\prod_{i=1}^{k}{\biggl(x-r(i)\biggr)^{m(i)}}\iota(x), \indent k=\tau(p),\indent \iota(x) \in \mathcal{I}_{\mathcal{R}}[x].$$ However, $q(x) \in \mathcal{V}(p(x))$ implies $q(x)\vert p(x)$ for any $q(x)\in \mathcal{V}(p(x))$ and also $\textbf{Deg}\biggl(q(x)\biggr)\ge \textbf{Deg}\biggl(k(x)\biggr)$ for any other $k(x)\vert p(x)$. Thus, one must have $q(x) = c\prod_{i=1}^{k}{\biggl(x-r(i)\biggr)^{m(i)}}$ or $q(x) = \prod_{i=1}^{k}{\biggl(x-r(i)\biggr)^{m(i)}}$. Since, $|q(x)|\ge |k(x)| $ for any other divisor $k(x)~\vert~ p(x)$ we must conclude that $$q(x)= c \prod_{i=1}^{k}{(x-r(i))^{m(i)}}=p(x)$$
\end{proof}
\end{lemma}

\begin{theorem}\label{Thm6}
let $\mathcal{R}$ be an integral domain with $\mathcal{R}\subset \mathbb{C}$. For every $p\in \mathcal{N}_{\mathcal{R}}[x]$ with degree $d= \textbf{Deg}(p)$ there exist a value $k\in \mathbb{N}$ and a set of functions $\{f_0,f_1,\cdots, f_d\}$ with $f_i:~\mathbb{C}^k\rightarrow \mathcal{R}$ for each $0\le i\le d$ such that

\begin{equation}\label{Eq22}
p(x) = \sum_{i=0}^{d}{f_i(\textbf{r})x^i}, \indent \textbf{r} = (r_1,r_2,...,r_k)
\end{equation}

with $\{r_1,r_2,\cdots r_k\}$ the unique roots to the polynomial and $k=\tau(p)$. 

\begin{proof}
The proof is inductive and a consequence of the structure of the ring $\mathcal{R}[x]$. First let $\mathcal{R}\subseteq \mathbb{C}$ and take $p\in \mathcal{N}_{\mathcal{R}}[x]\subset \mathbb{C}[x]$ to be a non-constant nicely factored polynomial. According to the Fundamental Theorem of Algebra $p$ has at least one root in $\mathbb{C}$. Set $d= \textbf{Deg}(p)$ and let $\{r_1,r_2,...,r_k\}\subset \mathbb{C}$ denote the set of unique roots with $k=\tau(p)$. Then $p$ can be written as

\begin{equation}\label{Eq23}
p(x) = c\prod_{i=1}^{k}{\biggl(x-r_i\biggr)^{m_i}},\indent c\in \mathcal{R}.
\end{equation}

First we check the base case $k=1$ is satisfied notice

\begin{equation}\label{Eq24}
\biggl(x-r_i\biggr)^{m_1} = \sum_{j=0}^{m_1}{ {{m_1}\choose {j}}}(-r_1)^{m_1-j}x^{j}
\end{equation}

define $f_i:\mathbb{C}^k\rightarrow \mathcal{R}$ using the rule 

\begin{equation}\label{Eq25}
f_i(\textbf{r})=  c{{m_1}\choose {j}}(-r_1)^{m_1-j}
\end{equation}

then we can write Eq.(\ref{Eq24}) as 

\begin{equation}\label{Eq26}
\sum_{j=0}^{m_1}{ {{m_1}\choose {j}}}(-r_1)^{m_i-j}x^{j} = \sum_{j=0}^{m_i}{f_j(\textbf{r})x^{j} }
\end{equation}

thus the base case is satisfied. Next assume the hypothesis is true for $1\le l\le k-1$ then there exist a set of function $\{g_0,g_1,\cdots, g_{k-1}\}$ with $g_i:\mathbb{C}^k\rightarrow \mathcal{R}$ for $1\le i\le k$ such that 

$$\prod_{i=1}^{k-1}{\biggl(x-r_i\biggr)^{m_i}}=\sum_{i=0}^{k-1}{g_i(\textbf{r})x^i}.$$

Define $f_i(\textbf{r})$ by the rule 

\begin{equation}\label{Eq27}
f_i(\textbf{r}) = c\sum_{\substack{i+j=n \\ 0 \leq i, j \leq m_k}}^{}{g_i(\textbf{r}){{m_i}\choose {j}}(-r_i)^{m_k-j}}
\end{equation}

and write 

\begin{equation}\label{Eq28}
\begin{split}
p(x) &  = c\prod_{i=1}^{k}{\biggl(x-r_i\biggr)^{m_i}}.\\
& = c\prod_{i=1}^{k-1}{\biggl(x-r_i\biggr)^{m_i}}\biggl(x-r_k\biggr)^{m_k}\\
& = c \biggl(\sum_{i=0}^{k-1}{g_i(\textbf{r})x^i}\biggr)\biggl(x-r_k\biggr)^{m_k}\\
& = c\biggl(\sum_{i=0}^{k-1}{g_i(\textbf{r})x^i}\biggr)\biggl(\sum_{j=0}^{m_k}{ {{m_k}\choose {j}}}(-r_i)^{m_k-j}x^{j} \biggr)\\
& = \sum_{i=0}^{k-1}{c\biggl(g_i(\textbf{r})\sum_{j=0}^{m_k}{ {{m_k}\choose {j}}}(-r_i)^{m_k-j}x^{i+j} \biggr)}\\
& = \sum_{n=0}^{k}{f_n(\textbf{r})x^i}\\
\end{split}
\end{equation}

By induction the theorem holds for all $k$.

\end{proof}
\end{theorem}

\begin{example}\label{Ex2}(5th degree polynomials roots) Let $p(x) \in \mathbb{Z}[x]$ with $\textbf{Deg}(p) = 5$, and suppose $p$ has $3$ unique roots $\{r_1, r_2, r_3\}\subset \mathbb{C}$. Then $p$ can be represented using its factorization

\begin{equation}\label{Eq29}
p(x) = c(x - r_1)(x - r_2)(x - r_3)\iota(x), \quad c \in \mathbb{Z}, \quad \iota(x) \in \mathcal{I}_{\mathbb{Z}}[x], \quad \textbf{Deg}(\iota) = 2.
\end{equation}

Set $\textbf{r} = (r_1,r_2,r_3)\in \mathbb{C}^3$ and define $g_1(\textbf{r}) = c$, $g_2(\textbf{r}) = -c(r_1+r_2+r_3)$, $g_3(\textbf{r}) = c(r_3r_1+r_3r_2-r_1r_2)$, and $ g_4(\textbf{r}) = cr_1r_2r_3$. Write the irreducible $i(x)$ using its coefficients $i(x) = a_1x^2+a_2x + a_3$ for $(a_1,a_2,a_3)\in \mathbb{Z}^3$ and consider the set of functions $\{f_0,f_1,f_2,f_3,f_4,f_5\}$ defined by the rules $$f_0(\textbf{r}) = a_1 g_1(\textbf{r}),$$ $$f_1(\textbf{r}) = a_1g_1(\textbf{r})+a_2g_2(\textbf{r}),$$ $$f_2(\textbf{r}) = a_1g_3(\textbf{r})+ a_2g_3(\textbf{r})+a_3g_1(\textbf{r}),$$ $$f_3(\textbf{r}) = a_1g_4(\textbf{r})+ a_2g_3(\textbf{r})+ a_3g_3(\textbf{r})),$$ $$f_4(\textbf{r}) = a_2g_4(\textbf{r})+a_3g_3(\textbf{r}),$$ $$f_5(\textbf{r}) = a_3g_4(\textbf{r}).$$ Then Eq(\ref{Eq29}) produces

\begin{equation}\label{Eq30}
\begin{split}
 p(x)  & = c(x - r_1)(x - r_2)(x - r_3)\iota(x)\\
& = \biggl(\textbf{g}_1(\textbf{r})x^3 + \textbf{g}_2(\textbf{r})x^2 + g_3(\textbf{r})x+ g_4(\textbf{r})\biggr)\iota(x)\\
 & = \biggl(\textbf{g}_1(\textbf{r})x^3 + \textbf{g}_2(\textbf{r})x^2 + g_3(\textbf{r})x+ g_4(\textbf{r})\biggr)\biggl(a_1x^2+a_2x + a_3\biggr)\\
 & = f_0(\textbf{r})x^5 + f_1(\textbf{r})x^4 + f_2(\textbf{r})x^3+ f_3(\textbf{r})x^2 + f_4(\textbf{r})x + f_5(\textbf{r}).
\end{split}
\end{equation}
\end{example}

\begin{example}\label{Ex3}(6th degree polynomials roots)
Let $p(x) \in \mathbb{Z}[x]$ with $\textbf{Deg}(p) = 6$, and suppose $p$ has $6$ unique roots $\{r_1, r_2, r_3, r_4, r_5, r_6\}$. Define $$g_1(\textbf{r})= -(r_1+r_2+r_3+r_4),$$ $$g_2(\textbf{r})= \biggl(r_3r_1+r_3r_2+ r_4r_1+r_4r_2+r_4r_3-r_1r_2\biggr),$$  $$g_3(\textbf{r}) =  -\biggl(r_4r_3r_1+r_4r_3r_2+ r_1r_2r_3-r_4r_1r_2\biggr),$$ and $$g_4(\textbf{r}) =-r_1r_2r_3r_4.$$ Next define the set of functions $\{f_0,f_1,f_2,f_3,f_4,f_5\}$  by the rules 

$$f_0(\textbf{r}) = c,$$ 

$$f_1(\textbf{r}) = c\biggl(g_1(\textbf{r})-r_5- r_6\biggr),$$

$$f_2(\textbf{r}) = c\biggl(g_2(\textbf{r})-r_5g_1(\textbf{r})-r_6g_1(\textbf{r})-r_6r_5)\biggr), $$ 
 
$$f_3(\textbf{r}) = c\biggl((g_3(\textbf{r})-r_5g_2(\textbf{r})-r_6g_2(\textbf{r})-r_6r_5g_1(\textbf{r})\biggr)), $$
 
$$f_4(\textbf{r}) = c\biggl((g_4(\textbf{r})-r_5g_3(\textbf{r}))- r_6(g_3(\textbf{r})-r_5g_2(\textbf{r}))\biggr), $$ 
  
$$f_5(\textbf{r}) =  -c\biggl(r_5g_4(\textbf{r}+r_6(g_4(\textbf{r})-r_5g_3(\textbf{r}))\biggr),$$ 
   
$$f_6(\textbf{r}) = cr_6r_5g_4(\textbf{r}).$$ 

Then $p$ can be represented using its factorization

\begin{equation}\label{Eq31}
\begin{split}
p(x) & = c(x - r_1)(x - r_2)(x - r_3)(x - r_4)(x - r_5)(x - r_6), \quad c \in \mathbb{Z}.\\
& = c\biggl(x^2-(r_1+r_2)x -r_1r_2\biggr) (x - r_3)(x - r_4)(x - r_5)(x - r_6)\\
& = c\biggl(x^3-(r_1+r_2+r_3)x^2+(r_3(r_1+r_2)-r_1r_2)x-r_1r_2r_3\biggr)(x-r_4)(x-r_5)(x-r_6)\\
& = c \biggl(x^4+ g_1(\textbf{r})x^3+ g_2(\textbf{r})x^2 + g_3(\textbf{r})x+g_4(\textbf{r})\biggr)(x-r_5)(x-r_6)\\
& = c \biggl(x^5 + (g_1(\textbf{r})-r_5)x^4+ (g_2(\textbf{r})-r_5g_1(\textbf{r}))x^3 + (g_3(\textbf{r})-r_5g_2(\textbf{r}))x^2\\
& ~~~~~~~~~~~~~~~~~~~~~~ + (g_4(\textbf{r})-r_5g_3(\textbf{r}))x-r_5g_4(\textbf{r})\biggr)(x-r_6)\\
& = f_0(\textbf{x})x^6 + f_1(\textbf{r})x^5 + f_2(\textbf{r})x^4 + f_3(\textbf{r})x^3 + f_4(\textbf{r})x^2 +f_5(\textbf{r})x + f_6(\textbf{r})
\end{split}
\end{equation}
\end{example}

\begin{example}\label{Ex4}
Determine the roots of the polynomial $p(x) = 2x^3-(\pi^2) x+ \dfrac{\pi}{2}$ in $\mathbb{R}$. We will separate the problem into three cases.\\

\textbf{Case I:} Assume $p$ has one root $\{r\}$ with $\Upsilon(p,r) = 3$ then write 

\begin{equation}\label{Eq32}
\begin{split}
p & = c(x-r)^3 \\
& = c(x-r)(x-r)(x-r) \\
& = c(x^2-2rx-r^2)(x-r)\\
& = cx^3 + c(-3r)x^2+c(2r-r^2)x+cr^3
\end{split}
\end{equation}

for some $c$. However, two polynomials equal only if there coefficients are equal and so $$cx^3 + c(-3r)x^2+c(2r-r^2)x+cr^3 = p(x) = 2x^3-(\pi^2) x+ \dfrac{\pi}{2}  $$

implies we must have 

\begin{subequations}\label{Eq33}
\begin{align}
c &= 2 \tag{33a} \\
c(-3r) &= 0 \tag{33b} \\
c(2r - r^2) &= -\pi^2 \tag{33c} \\
cr^3 &= \frac{\pi}{2} \tag{33d}
\end{align}
\end{subequations}

Notice Eq.(\ref{Eq33}.a) and Eq.(\ref{Eq33}.b) implies $r=0$ which is not possible. We conclude that $p$ cannot have one root. \\

\textbf{Case II:} Assume $p$ has two unique roots $\{r_1,r_2\}$ with $\Upsilon(p,r_1)=2$ or  $\Upsilon(p,r_2)=2$ then write 

\begin{equation}\label{Eq34}
p(x) = c(x-r_1)(x-r_2)^2,\indent \text{or}\indent p(x)= c(x-r_1)^2(x-r_2)
\end{equation}

for some constant $c\in \mathbb{R}$ with $c\neq 0$. First assume $p$ satisfies the first. Then 

\begin{equation}\label{Eq35}
\begin{split}
p(x) & = c(x-r_1)(x-r_2)^2\\
& = c(x-r_1)(x^2-2r_2x-r_2^2)\\
& = cx^3 + c(-2r_2-r_1)x^2 + c(-r_2^2+2r_2)x+cr_1r_2^2
\end{split}
\end{equation}

Consequently, we must have 

\begin{subequations}\label{Eq36}
\begin{align}
c &= 2 \tag{34a} \\
c(-2r_2-r_1)&= 0 \tag{34b} \\
c(-r_2^2+2r_2) &= -\pi^2 \tag{34c} \\
cr_1r_2^2 &= \frac{\pi}{2} \tag{34d}
\end{align}
\end{subequations}

Eq.(\ref{Eq34}b) implies $r_1 = -2r_2$. Substituting this result into Eq.(\ref{Eq34}d) yields $(-2r_2)r_2^2 = \frac{\pi}{2} $ or  $r_2^3 = -\dfrac{\pi}{4} \Rightarrow r_2 =  -\sqrt[3]{\dfrac{\pi}{4}}$. Using this result in Eq.(\ref{Eq34}c) gives ones 

$$\biggl((\sqrt[3]{\dfrac{\pi}{4}})^2- 2\sqrt[3]{\dfrac{\pi}{4}}\biggr) = \dfrac{\pi^2}{2}$$

which is false. The same situation applies is $\Upsilon(p,r_2)=2$ and $\Upsilon(p,r_1) = 1$ as it will not affect the results. Therefore we conclude that $p$ has more then two roots or no roots at all as case I and case II are both false. \\

\textbf{Case III:} Assume $p$ has three unique roots namely $\{r_1,r_2,r_3\}$ and write 

\begin{equation}\label{Eq37}
\begin{split}
p(x) & = c(x-r_1)(x-r_2)(x-r_3)\\
& = c(x^2-(r_1+r_2)x+r_1r_2)(x-r_3)\\
& = cx^3 -c(r_1+r_2+r_3)x^2 + c\biggl(r_1r_2+r_3(r_1+r_2)\biggr)x -cr_1r_2r_3
\end{split}
\end{equation}

As a result we must have

\begin{subequations}\label{Eq38}
\begin{align}
c &= 2 \tag{36a} \\
-c(r_1+r_2+r_3) &= 0 \tag{36b} \\
  c\biggl(r_1r_2+r_3r_1+r_3r_2\biggr) &= -\pi^2 \tag{36c} \\
 -cr_1r_2r_3 &= \frac{\pi}{2} \tag{36d}
\end{align}
\end{subequations}

Notice Eq.(\ref{Eq36}b) gives $r_1+r_2+r_3=0$ and so $r_1 = -(r_2+r_3)$. From Eq.(\ref{Eq36}d) we have 

\begin{equation}\label{Eq39}
2(r_2^2r_3+r_2r_3^2)  = \dfrac{\pi}{2}.
\end{equation}

Similarly, substituting $r_1 = -(r_2+r_3)$ into Eq.(\ref{Eq36}c) produces

\begin{equation}\label{Eq40}
\begin{split}
2\biggl(-(r_2+r_3)r_2-r_3(r_2+r_3)+r_3r_2\biggr)   & = -\pi^2\\
& \Downarrow\\
-2\biggl(-(r_2+r_3)^2+r_3r_2\biggr) & =  \pi^2\\
& \Downarrow\\
\dfrac{\sqrt{-2\biggl(-(r_2+r_3)^2+r_3r_2\biggr)}}{2} & =  \dfrac{\pi}{2}
\end{split}
\end{equation}

Comparing the results of Eq.(\ref{Eq39}) and Eq.(\ref{Eq40}) yields 

\begin{equation}\label{Eq41}
\begin{split}
2(r_2^2r_3+r_2r_3^2) & = \dfrac{\sqrt{-2\biggl(-(r_2+r_3)^2+r_3r_2\biggr)}}{2}\\
& \Downarrow\\
8(r_2^2r_3+r_2r_3^2)^2 -(r_2+r_3)^2+r_3r_2 & = 0\\
& \Downarrow\\
8(r_2r_3 (r_2+r_3))^2 -(r_2+r_3)^2+r_3r_2 & = 0\\
& \Downarrow\\
8 (\epsilon^2 \kappa^2) - \kappa^2+\epsilon &= 0,\indent \kappa  : = r_2+r_3, ~~\epsilon : = r_2r_3.
\end{split}
\end{equation}

Notice,

\begin{equation} \label{Eq42}
\begin{split}
8 (\epsilon^2 \kappa^2) - \kappa^2+\epsilon & = 0 \\
\Downarrow\\
\kappa^2 (8\epsilon^2-1)+\epsilon & = 0 \\
\Downarrow\\
\kappa & = \pm \sqrt{\dfrac{-\epsilon}{8\epsilon^2-1}}
\end{split}
\end{equation}

The equation $\kappa = \pm \sqrt{\dfrac{-\epsilon}{8\epsilon^2-1}} $ gives one a relationship between the sum $r_2+r_3$ and the product $r_2r_3$. Substituting either case back into the first part of Eq.(\ref{Eq41}) yields 

\begin{equation}\label{Eq43}
\begin{split}
2\epsilon \sqrt{\dfrac{-\epsilon}{8\epsilon^2-1}} & = \dfrac{\sqrt{-2 \biggl(\dfrac{-\epsilon}{8\epsilon^2-1}+\epsilon \biggr)}}{2} \\
& \Downarrow\\
\dfrac{8\epsilon^2}{8\epsilon^2-1}  & =  \biggl(\dfrac{\epsilon(2-8\epsilon^2)}{1-8\epsilon^2}\biggr) \\
& \Downarrow\\
8\epsilon^2 (1-8\epsilon^2)  & = \epsilon(8\epsilon^2-1)(2-8\epsilon^2)\\
& \Downarrow\\
8\epsilon^2 -64\epsilon^4  & = (8  \epsilon^2 - \epsilon)(2-8\epsilon^2)\\
& \Downarrow\\
8\epsilon^2 -64\epsilon^4  & = 16\epsilon^2-64\epsilon^4 - 2\epsilon+ 8\epsilon^3)\\
& \Downarrow\\
0 & =  \epsilon(8\epsilon^2+8\epsilon - 2)\\
\end{split}
\end{equation}

Thus $\epsilon =0$ or $(8\epsilon^2+8\epsilon - 2)=0$. We note $\epsilon =0 \rightarrow r_2r_3=0$ since we are working in a integral domain it has no zero divisors this gives $r_2=0$ or $r_3=0$ which is not possible. Thus we must have $(8\epsilon^2+8\epsilon - 2)=0$ which implies 

\begin{equation}\label{Eq44}
\epsilon = \dfrac{-8+ \sqrt{128}}{16},\indent \text{or} \indent \epsilon = \dfrac{-8- \sqrt{128}}{16}
\end{equation}

however, $\kappa\in \mathbb{R}$ (as $\kappa=r_1+r_2$) implies $\epsilon<0$ and so we must have $\epsilon = \dfrac{-8- \sqrt{128}}{16}$ and consequently

\begin{equation}\label{Eq45}
 r_2 = \dfrac{1}{r_3}\biggl(\dfrac{-8- \sqrt{128}}{16}\biggr).
\end{equation}

Substituting the result of Eq.(\ref{Eq45}) into Eq.(\ref{Eq39}) gives an explicit result $r_3$:

\begin{equation}\label{Eq46}
\begin{split}
2(r_2r_3 (r_2+r_3)) & =  \dfrac{\pi}{2}\\
& \Downarrow\\
2(\dfrac{1}{r_3}\biggl( \dfrac{-8+ \sqrt{128}}{16}\biggr)(\biggl( \dfrac{-8+ \sqrt{128}}{16}\biggr)+r_3^2)) & =  \dfrac{\pi}{2}\\
&\Downarrow\\
(\biggl( \dfrac{-8- \sqrt{128}}{16}\biggr)(\biggl( \dfrac{-8- \sqrt{128}}{16}\biggr)+r_3^2)) -   r_3\dfrac{\pi}{4} & = 0\\
\end{split}
\end{equation}

The last equation is a quadratic equation in $r_3$. We note for $a\in \mathbb{R}$ the roots of $a(a+r_3)-r_3\dfrac{\pi}{4} = a^2 + ar_3^2 -r_3 \dfrac{\pi}{4}$ are provided by $r_3 = \dfrac{\dfrac{\pi}{4} \pm \sqrt{\dfrac{\pi^2}{16}-4a^3}}{2a}.$ With respect to Eq.(\ref{Eq46}) $a=\dfrac{-8- \sqrt{128}}{16}$ we must have $\dfrac{\pi^2}{16}-4a^3\ge 0$. 
Note $r_3\ge 0$ implies $r_2\le 0$ from Eq.( \ref{Eq45}) but recall $r_1 = -(r_2+r_3)$ so $r_2\le 0$ whenever $r_3\ge 0$ . So the sign of $r_2$ and $r_3$ are different. In order to fully determine the correct value of $r_3$ we ask the question 
if $a= \dfrac{-8- \sqrt{128}}{16} $ then $\dfrac{\dfrac{\pi}{4} + \sqrt{\dfrac{\pi^2}{16}-4a^3}}{2a}\ge 0$ or $\dfrac{\dfrac{\pi}{4} -\sqrt{\dfrac{\pi^2}{16}-4a^3}}{2a}\ge 0$ but not both. Therefore, the positive root is $r_3 = \dfrac{\dfrac{\pi}{4} - \sqrt{\dfrac{\pi^2}{16}-4a^3}}{2a}$. Thus, if the polynomial $p(x) = 2x^3-(\pi^2) x+ \dfrac{\pi}{2}$ has any roots at all it must have three provided by 

\begin{equation}\label{Eq47}
\begin{split}
r_3   & =\dfrac{\dfrac{\pi}{4} - \sqrt{\dfrac{\pi^2}{16}-4a^3}}{2a}, \indent a = \dfrac{-8- \sqrt{128}}{16} \\
r_2 & = \biggl(\dfrac{-8- \sqrt{128}}{16}\biggr)/\biggl(\dfrac{\dfrac{\pi}{4} - \sqrt{\dfrac{\pi^2}{16}-4a^3}}{2a}\biggr)\\
r_1 & = -\biggl(r_2+r_3\biggr)
\end{split}
\end{equation}

It is not difficult to see $p(-\pi) = -\pi^3 + \dfrac{\pi}{2} <0$ and $p(\pi) = \pi^3 + \dfrac{\pi}{2}>0$ so $p$ has at least one root as it is continuous on $\mathbb{R}$. From the combination of case (I)- case (III) we conclude that Eq.(\ref{Eq47}) prescribes them analytically and there are exactly three.  
\end{example}


\begin{theorem}\label{Thm7}
For $p\in \mathcal{R}[x]\subseteq \mathbb{C}[x]$ there is set of solvable algebraic systems that fully determine the roots of $p$ in $\mathbb{C}$.

\begin{proof}
Let $p\in \mathcal{R}[x]\subseteq\mathbb{C}[x]$ and write $p= \sum_{i=1}^{d}{a_ix^i}$ for some function $a_i:\mathbb{N}\rightarrow \mathbb{N}$ and set $\tau(p) = k$. Corollary to theorem \ref{Thm5} there exist a set a functions $\biggl\{f_0,f_1,\cdots,f_d\biggr\}$ with $f_i:\mathbb{C}^k\rightarrow \mathcal{R}$ for $0\le i\le d$ such that

\begin{equation}\label{Eq48}
\begin{bmatrix}
f_0(\textbf{r})\\
f_1(\textbf{r}) \\
\vdots\\
f_d(\textbf{r}) 
\end{bmatrix} = \begin{bmatrix}
a_0\\
a_1 \\
\vdots \\
a_d .
\end{bmatrix}
\end{equation}

This is a vector valued equation over $\underbrace{\mathcal{R}\times \mathcal{R}\times \cdots \times \mathcal{R}}_{d-\text{times}}$ which is solvable for a unique value $\textbf{r}\in \mathbb{C}^k$ as $k\le d$.
\end{proof}
\end{theorem}

\subsubsection{Time Complexity}
A key challenge in our work lies in determining whether finding a root of a polynomial can be as computationally efficient as verifying one. Verifying that a value \( a \) satisfies \( p(a) = 0 \) is straightforward and typically (not not always) involves evaluating \( p(a) \) directly using the polynomial's definition. In contrast, finding the entire set of roots for a given polynomial \( p(x) \) is considerably more complex, particularly for high-degree polynomials, as no general solution is known for arbitrary cases. Even if such a solution existed, computational efficiency would remain a concern. Polynomial evaluation can be significantly accelerated using \textit{Horner’s method} \cite{Burrus2003, Knuth1962, Cajori1911}, an efficient algorithm that reduces the computational cost. However, finding polynomial roots requires more sophisticated approaches, and numerous algorithms have been developed for this purpose, each with distinct strengths and applications:

\begin{itemize}
    \item \textbf{Newton's Iterative Method}: Known for its rapid convergence, this method has been successfully applied to polynomials with degrees exceeding a million \cite{Schleicher2017}.
    \item \textbf{Francis QR Algorithm}: A robust approach for eigenvalue problems, with generalizations extending its applicability \cite{Parlett1968, Watkins2011, Bunse-Gerstner1985}.
    \item \textbf{Aberth’s Method}: A globally convergent iterative method for locating all polynomial roots simultaneously \cite{Aberth1998, Aberth2007}.
    \item \textbf{Power Method} and \textbf{Bernoulli’s Method}: Effective in specific contexts, offering unique advantages \cite{McNamee2013, Moller}.
    \item \textbf{Jenkins–Traub Algorithm}: Renowned for its accuracy and reliability, with several generalizations enhancing its scope \cite{Jenkins, Jenkins2, Traub, Ford1977, Ellis}.
\end{itemize}

Recent advancements include the application of matrix factorization techniques, such as QR decomposition \cite{Fry1945, Pan2011, Bini2005}, which have accelerated root-finding for large polynomials. In the context of the \(\textbf{P}~\text{vs.}~\textbf{NP}\) problem, the computational complexity of determining polynomial roots is encapsulated in the following question:\\

\textbf{Problem 2}: (Roots of Polynomials)  
For \( p(x) \in \mathbb{R}[x] \), can one determine the set \( \{a \in \mathbb{C} : p(a) = 0\} \) in polynomial time?\\


Regardless of a polynomial's degree, the minimal time required to evaluate a root and verify its validity is fundamentally determined by the minimal time it takes to evaluate the polynomial $p$ at any point in its domain (some points may be easier then others i.e $1$, $0$, or integer powers $2^z$, $3^z$, ... etc ) which can be done using circuit complexity theoretical approaches \cite{Kabanets}. However, as the degree of the polynomial tends to infinity, it is reasonable to suspect that determining its roots becomes increasingly challenging for any algorithm similarity to how the complexity of the sliding puzzle appears to increase as the order or dimension of a tensor sliding puzzle increases. We will formalize this concept in the next section, but for now, it is evident that addressing these issues provides valuable insight into the problem at hand. Formally, we will need to discuss the minimal computable function $\text{COMP}_m(f)$ of a general solution for any deterministic or non-deterministic algorithm and relate it to the minimal computable function of the versifier function which will usually involve parameters in the statement considered using circuit complexity techniques combines with set theoretical axioms and to determine if $\text{COMP}_m(f)$ can be effectively mapped to in polynomial time.

\section{Generalization to $\textbf{P}$ vs \textbf{NP}}\label{sec3}
In this section we extend the concepts and problems  introduced in Sections \ref{Sec1} and \ref{Sec2} to explore their implications for the $\textbf{P}$ vs $\textbf{NP}$ problem, in particular showing that verifying sliding puzzles can be done in polynomial time. Cook's seminal work \cite{cooks3}, remains one of the most significant open questions in theoretical computer science. While there is some debate regarding the exact origins of the problem, it has inspired a diverse array of approaches, including thermodynamic models \cite{Neukart1, Neukart2}, algebraic geometry methods for NP-hard problems \cite{Cheng2008, Borodin1983}, geometric complexity theory (GCT) \cite{Mulmuley1}, quantum computing \cite{Aaronson, Chapman2022}, and circuit theory \cite{Circuit0,Circuit1,Circuit2,Circuit3,Circuit4}  \cite{Circuit5}. In Computability theory the notion of a computable function is now standard \cite{Fenstad}. For our purpose a Turing machine $\mathcal{M}$, first defined by Alan Turing \cite{Turing1936}, serves as a theoretical framework for modeling computation. It processes input strings from a finite alphabet $\mathcal{B} = \{b_1, b_2, \dots, b_k\}$ of length $k$, where $\mathcal{L}$ represents the set of all possible strings formed from $\mathcal{B}$

\begin{equation}\label{Eq49}
\mathcal{L} = \biggl\{b_{\rho(1)} b_{\rho(2)}b_{\rho(3)} \cdots b_{\rho(l)}~\vert~~\rho:\mathbb{N}\rightarrow \{1,2,\cdots,k\},~ 1\le l\le k  \biggr\}.
\end{equation}

The Turing machine $\mathcal{M}$ computes a time-consuming function $\mathcal{K} : \mathcal{L}  \to \mathbb{R} \cup \{\infty\}$, where $\mathcal{K}(\sigma) = t$ for input $\sigma \in \mathcal{L} $. If $t \in \mathbb{R}$, we say $\mathcal{M}$ accepts $\sigma$, and if $t = \infty$ we say it rejects $\sigma$. Allow $\mathcal{M}^{\star}$ to denote the set of all Turing machines

\begin{equation}\label{Eq50}
\mathcal{M}^{\star}=\bigcup_{\text{$\mathcal{M} $ s.t $\mathcal{M}$ is a Turing model}}^{}{\mathcal{M}}
\end{equation}

Each Turing machine $\mathcal{M} \in \mathcal{M}^\star$ is associated with a decision function $\mathcal{J} : \mathcal{L}  ~\times ~\mathcal{M}^\star \to \{1, 0\}$, defined as:

\begin{equation}\label{Eq51}
\mathcal{J}(\sigma,\mathcal{M}) = \begin{cases}
1 & \text{if $\mathcal{M}$ accepts $\sigma$}\\
0 & \text{if $\mathcal{M}$ rejects $\sigma$}.
\end{cases}
\end{equation}

\begin{definition}\label{Def14}
A \textit{computable function} is an element $f\in \mathcal{L} $ such that there exist a Turing model $\mathcal{M}\in \mathcal{M}^\star$ so that $\mathcal{J}(f,\mathcal{M}) =1$. It is simply element of a base that is accepted by some Turing model. These are also all algorithm we can compute in a finite amount of time.
\end{definition}

\begin{definition}\label{Def15}
The \textit{time complexity} of a Turing machine \( M \) and a computable function \( f \) is denoted by \( \mathcal{T}_M(f) \), and represents the number of steps \( M \) takes to compute \( f(x) \) for a given input \( x \). That is the number of decision processes needed to compute $f$.
\end{definition}

\begin{definition}\label{Def16}
Let \( f \) be a computable function and \( \mathcal{M} \) be a Turing machine the function $\textbf{Res}_{\mathcal{M}}(f)$ denotes the result of the function (which is a element of the set $\mathcal{L} $).
\end{definition}

\begin{definition}\label{Def17}
Let \( f \) be a computable function and \( \mathcal{M} \) be a Turing machine. The set of \textit{ minimal computable} functions with respect to $f$ is provided by the map

\begin{equation}\label{Eq52}
(f,\mathcal{M}) \mapsto \{ g \in \mathcal{L}  ~: ~\textbf{Res}_{\mathcal{M}}(g) =  \textbf{Res}_{\mathcal{M}}(f), ~~\wedge~~ \mathcal{T}_{\mathcal{M}}(g) \leq \mathcal{T}_{\mathcal{M}}(f) \}
\end{equation}

We denote the image as $\text{COMP}_{\mathcal{M}}(f) $. This is simply the set of all computable function that can be computed in the same amount or less time then any other function accepted by a Turing machine $\mathcal{M}$ that also accepts $f$.
\end{definition}

\subsubsection{Sliding tiles and $\textbf{P}$ vs \textbf{NP}}
In view of the sliding puzzle problems, as the dimension of different puzzles tends to be infinitely large we suspect the size of solution space becomes exhaustively large causing the minimal amount of steps for any algorithm or computable function to compute the "general solution" to be exponential relative to the minimal of steps to compute any finite sequence of moves to any puzzle (as we only need to swap one element in a tensor with the empty set which can be done in polynomial time by searching ever possible element of the tensor and verifying if its the empty set and if so swapping it with the correct tile producing a result). Verifying a seqeuence of moves produces the the goal state takes time compute each element in the sequence of moves but also at the last stage checking if we have reached the goal state.

\begin{lemma}\label{lemma4}
\begin{equation}\label{Eq53}
\mathcal{T}_{\mathcal{M}}(\rho_k(i,j))\le 4 \indent k\in \{1,2,3,4\}.
\end{equation}
\begin{proof}
$\rho_j$ is a case statement that takes at most four decisions to decide its result (the vector it returns). 
\end{proof}
\end{lemma}

\begin{lemma}\label{lemma5}
\begin{equation}\label{Eq54}
\mathcal{T}_{\mathcal{M}}(\tau_k(i,j))\le 17 \indent k\in \{1,2,3,4\}
\end{equation}

\begin{proof}
There are a maximum of four decisions that it takes to compute $\rho_k$ from Lemma \ref{lemma4}, but also $ \mathcal{T}_{\mathcal{M}}(\tau_k))\le 1+\mathcal{T}_{\mathcal{M}}(\lambda_k)$ and it takes at most four decisions to decide if $\tau_k$ maps to $(i,j)$ or $\lambda_k(i,j)$ because of the case statement Eq.(\ref{Equation 6}). Lastly, it takes a maximum of four decisions to decide what vector $\lambda_k(i,j)$  actually is. But $\mathcal{T}_{\mathcal{M}}(\lambda_k(i,j))\le 4$ so it cant take more $16$ decisions (the product $4*4$ four decision max for each case statement) to determine $\lambda_k$ output and then $1+16$ decisions to decide the result $\tau_k$ meaning $\mathcal{T}_{\mathcal{M}}(\tau_k(i,j))\le 17$.

\end{proof}

\end{lemma}

\begin{lemma}\label{lemma6}
\begin{equation}\label{Eq55}
\mathcal{T}_{\mathcal{M}}(\phi_{j}(\textbf{A}))\le 5 + \mathcal{T}_{\mathcal{M}}(\Sigma_j(\textbf{A})) ,\indent j\in\{1,2,3,4\}.
\end{equation}

\begin{proof}
Any tile movement $\phi_j$ takes a maximum of five decisions before applying $\Sigma_j$ and $\Sigma_j$ must be applied after.
\end{proof}
\end{lemma}

\begin{lemma}\label{lemma7}
\begin{equation}\label{Eq56}
\mathcal{T}_{\mathcal{M}}(\Sigma_j(\textbf{A}))\le 1+ \mathcal{T}_{\mathcal{M}}(\rho_j(\textbf{A}))+\mathcal{T}_{\mathcal{M}}(\tau_j(\textbf{A})) ,\indent j\in\{1,2,3,4\}.
\end{equation}

\begin{lemma}\label{lemma8}
The general solution $f$ to all $n\times n$ sliding tiles in Definition \ref{Def11} is computable.
\begin{proof}

The proof follows as a corollary of Theorem \ref{Thm3} or  Theorem \ref{Thm4} as we can exhaust the solution space. We use the first because it holds for more values of $n$ not because it is the best bound we found. We have 

\begin{equation}\label{Eq57}
|\mathcal{S}_n(\mathbb{N}\cup \{\varnothing\})|\le 4 \biggl(3^{\lfloor \dfrac{\log_4(n^2!)-1}{2}\rfloor }\biggr)\biggl(4^{\lfloor \dfrac{ \log_4(n^2!)-1}{2}}\rfloor\biggr)+4.
\end{equation}

Set $k=4 \biggl(3^{\lfloor \dfrac{\log_4(n^2!)-1}{2}\rfloor }\biggr)\biggl(4^{\lfloor \dfrac{ \log_4(n^2!)-1}{2}}\rfloor\biggr)+4$ for any given $n>1$. Consider the set of all possible consecutive application of the four possible legal transformations with the number of compositions less then $k$

\begin{equation}\label{Eq58}
\mathcal{P}^k = \biggl\{(\epsilon_1,\epsilon_2,\cdots, \epsilon_q)\mapsto  \biggl(\Phi_{\epsilon_1}\circ \Phi_{\epsilon_2}\circ \cdots \circ \Phi_{\epsilon_q}\biggr)~: ~~1\le q\le k , ~~\epsilon_i\neq \epsilon_j~~ \forall ~1\le j\le n , j\neq i \biggr\}.
\end{equation}

Since $f$ is the general solution it must map to the goal state. Since $\mathcal{P}^k$ contains all possible maps that will reach the goal state we know that $f\in \mathcal{P}^k$ but

\begin{equation}\label{Eq59}
|\mathcal{P}^k| = \sum_{l=1}^{k}{4^l} \le  4^{k+1} <\infty.
\end{equation}

Take a bijection $\mathcal{Z}:\{1,2,..., |\mathcal{P}|^k\} \rightarrow \mathcal{P}^k$ and consider the map

\begin{equation}\label{Eq60}
\mathcal{F}(i) = \begin{cases}
1 & \text{if $\mathcal{Z}(i)\circ \textbf{A} = \textbf{A}_{goal} $}\\
0 & \text{Otherwise} \\
\end{cases}
\end{equation}

and take the element $\sigma \in \mathcal{L}$ (combined with a Turing $\mathcal{M}$ that accepts $\sigma$) defined as 

\begin{equation}\label{Eq61}
\sigma = "\text{For each } 1\le i \le 4^k, \text{ if } \mathcal{F}(i) = 1, \text{ return } \mathcal{Z}(i), \text{ else continue.}"
\end{equation}

where of course we initialize the function $\mathcal{F}(i)$. This element $\sigma\in \mathcal{L}$ will terminate in at most $4^k$ steps but along each step we must compute $\mathcal{F}(i)$ for each $1\le i\le k$ as we are exhausting the solution space. Let us prove it more thoughtfully. Allow $g:M_n(\mathbb{N}\cup\{\varnothing\}) \times M_n(\mathbb{N}\cup\{\varnothing\})\rightarrow  \{0,1\}$ to be the computable function defined as

\begin{equation}\label{Eq62}
(\textbf{A}, \textbf{B})\mapsto \begin{cases}
1 & \text{if $\textbf{A}=\textbf{B}$}\\
0 & \text{Otherwise}
\end{cases}
\end{equation}

The total time it takes to compute $\sigma$ is

\begin{equation}\label{Eq63}
\begin{split}
\mathcal{T}_{\mathcal{M}}(\sigma)& \le \sum_{i=1}^{4^k}{\mathcal{T}_{\mathcal{M}}(\mathcal{F}(i))}\\
 & \le \sum_{i=1}^{4^k}{\biggl[ 1+ \mathcal{T}_{\mathcal{M}}(\mathcal{Z}(i)\circ \textbf{A})+ g(\mathcal{Z}(i)\circ \textbf{A}), \textbf{A}_{goal})\biggr]} \\
& \le 4^k + \sum_{i=1}^{4^k}{\textbf{Max}\biggl(\{\mathcal{T}_{\mathcal{M}}(f):~~ f\in \mathcal{P}^k\}\bigg)}+\sum_{i=1}^{4^k}{\mathcal{T}_{\mathcal{M}}(g(\mathcal{Z}(i)\circ \textbf{A}), \textbf{A}_{goal}))} \\
& \le 4^k + \sum_{i=1}^{4^k}{\textbf{Max}\biggl(\{\mathcal{T}_{\mathcal{M}}(f):~~ f\in \mathcal{P}^k\}\bigg)}+\sum_{i=1}^{4^k}{n^2} \\
& \le 4^k(1+n^2) + \sum_{i=1}^{4^k}{\textbf{Max}\biggl(\{\mathcal{T}_{\mathcal{M}}(f):~~ f\in \mathcal{P}^k\}\bigg)}\\
& \le 4^k(1+n^2)+ \sum_{i=1}^{4^k}{\mathcal{T}_{\mathcal{M}}\biggl((\Phi_{\epsilon_1}\circ \Phi_{\epsilon_2}\circ \cdots \circ \Phi_{\epsilon_k})\circ \textbf{A}\biggr)}\\
& \le 4^k(1+n^2)+ \sum_{i=1}^{4^k}{\sum_{j=1}^{k}{\mathcal{T}_{\mathcal{M}}(\phi_{\epsilon_j}\circ \textbf{A})}}\\
& \le 4^k(1+n^2)+ \sum_{i=1}^{4^k}{\sum_{j=1}^{k}{5 + \mathcal{T}_{\mathcal{M}}(\Sigma_j(\textbf{A})) }}\\
& \le 4^k(1+n^2)+ \sum_{i=1}^{4^k}{\sum_{j=1}^{k}{5 + 1+ \mathcal{T}_{\mathcal{M}}(\rho_j(\textbf{A}))+\mathcal{T}_{\mathcal{M}}(\tau_j(\textbf{A}))}} \\
& \le 4^k(1+n^2)+ \sum_{i=1}^{4^k}{\sum_{j=1}^{k}{5 + 1+ 4+17}} \\
& =  4^k(n^2+2) +27k
\end{split}
\end{equation}

We have shown $\mathcal{T}_{\mathcal{M}}(\sigma)\le 4^k(n^2+2) +27k$ which mean $f$ can be computed within $4^k(n^2+2) +27k$ many steps.
\end{proof}
\end{lemma}

\begin{definition}\label{Def18}
The function $\textbf{length}:~\mathcal{L}  \rightarrow \mathbb{N}$ takes in a string $\sigma \in \mathcal{L} $ and returns the length of the string $\sigma$. For example

\begin{equation}\label{Eq64}
\textbf{Length}(1!+456j)=7,\indent{and}\indent \textbf{Length}\biggl(y=3;~for(i=1;~i<k;~i++)\{~y=y+i;\}\biggr)=28.
\end{equation}

We do not include empty spaces as a string element.
\end{definition}

\begin{definition}\label{Def19}
A Computable function is said to run in polynomial-time if there exist a $p\in \mathbb{R}[x]$ such that that 

\begin{equation}\label{Eq65}
\mathcal{T}_{\mathcal{M}}(f) \le p(\textbf{Length}(f)).
\end{equation}
\end{definition}
\end{lemma}

\begin{theorem}(Verifying sliding tiles is in $\textbf{P}$)\label{Thm8}
Take $(n,k)\in \mathbb{N}\times \mathbb{N}$, $\textbf{A}\in \mathcal{S}_n(\mathbb{N}\cup \{\varnothing\})$, $\epsilon:\mathbb{N}\rightarrow \mathbb{N}$ to be a computable function, set $f$ to be the general solution to the sliding tiles as in Definition \ref{Def6}, and $h:\mathcal{S}_n(\mathbb{N}\cup \{\varnothing\})\times \mathbb{N}\rightarrow \mathcal{S}_n(\mathbb{N}\cup \{\varnothing\})$ to be the computable function defined by the rule 

\begin{equation}\label{Eq66}
(\textbf{A},k)\mapsto (\phi_{\epsilon(1)}\circ \phi_{\epsilon(2)}\circ \cdots \phi_{\epsilon(k)} )\circ \textbf{A}.
\end{equation}

Next allow $g:M_n(\mathbb{N}\cup\{\varnothing\}) \times M_n(\mathbb{N}\cup\{\varnothing\})\rightarrow  \{0,1\}$ to be the computable function defined as

\begin{equation}\label{Eq67}
(\textbf{A}, \textbf{B})\mapsto \begin{cases}
1 & \text{if $\textbf{A}=\textbf{B}$}\\
0 & \text{Otherwise}
\end{cases}
\end{equation}

and take $\gamma \in \text{Comp}_{\mathcal{M}}(f)$ then 
 
\begin{equation}\label{Eq68}
\mathcal{T}_\mathcal{M}(\gamma)\le\biggl(\textbf{length}(\gamma)^2+1\biggr)^{n^2+27k + 1}.
\end{equation}

\begin{proof}
Consider the three computable functions $\alpha\in \text{Comp}_{\mathcal{M}}(h) $, $\beta \in \text{Comp}_{\mathcal{M}}(g) $, and $\gamma \in \text{Comp}_{\mathcal{M}}(f)$. From Lemma \ref{lemma4}-\ref{lemma8} we see that

\begin{equation}\label{Eq69}
\begin{split}
\mathcal{T}_\mathcal{M}(\gamma) &   \le \mathcal{T}_\mathcal{M}\biggl(\alpha\biggr)  +\mathcal{T}_\mathcal{M}\biggl(\beta\biggr)  \\
& \le \sum_{j=1}^{k}{\mathcal{T}_\mathcal{M}(\phi_{\epsilon(j)}(\textbf{A}))}+\mathcal{T}_\mathcal{M}\biggl(g(h(\textbf{A}), \textbf{A}_{goal})\biggr)+1\\
& \le \sum_{j=1}^{k}{\biggl(5+\mathcal{T}_\mathcal{M}(\Sigma_j)\biggr)}+\mathcal{T}_\mathcal{M}\biggl(g(h(\textbf{A}), \textbf{A}_{goal})\biggr)+1\\
& \le \sum_{j=1}^{k}{\biggl(5+\biggl(1+\mathcal{T}_\mathcal{M}(\rho_j)+\mathcal{T}_\mathcal{M}(\tau_j)\biggr)\biggr)}+\mathcal{T}_\mathcal{M}\biggl(g(h(\textbf{A}), \textbf{A}_{goal})\biggr)+1\\
& \le \sum_{j=1}^{k}{\biggl(5+\biggl(1+4+17\biggr)\biggr)}+\mathcal{T}_\mathcal{M}\biggl(g(h(\textbf{A}), \textbf{A}_{goal})\biggr)+1\\
& = \sum_{j=1}^{k}{27}+\mathcal{T}_\mathcal{M}\biggl(g(h(\textbf{A}), \textbf{A}_{goal})\biggr)+1\\
& = n^2+27k + 1\\
& \le \biggl(x^2+1\biggr)^{n^2+27k + 1}\indent x\ge 1
\end{split}
\end{equation}

set $x= \textbf{length}(\gamma)\ge 1$ then we have $\mathcal{T}_{\mathcal{M}}(f) \le \biggl(\textbf{length}(\gamma)^2+1\biggr)^{n^2+27k + 1}.$ Since $\mathbb{R}[x]$ is a ring and closed under multiplication we know that

\begin{equation}\label{Eq70}
\biggl(\textbf{length}(\gamma)^2+1\biggr)^{n^2+27k + 1} = \prod_{i=1}^{n^2+27k + 1}{\biggl(\textbf{length}(\gamma)^2+1\biggr)}
\end{equation}

is a polynomial in $\gamma$. Thus, verifying a solution to a sliding tile can be accomplished in polynomial time.
\end{proof}
\end{theorem}

\subsubsection{Roots, algebraic system, and $\textbf{P}$ vs \textbf{NP}}
So far, we have discussed some of the time complexity issues related to the sliding puzzle problem but we have not discussed the time complexity in relation to the roots of polynomials. In this short section, we will explore in detail the problem concerning time complexity of finding versus verifying roots in an arbitrary ring by generalizing the algebraic system in Theorem \ref{Thm6} and \ref{Thm7}. This means examining the time that it takes to set up the algebraic system and then, secondly, the time it takes to solve it and get the output. Nonetheless, theoretical methods may be able to show that a "general solution" to the algebraic system Eq.(\ref{Eq48}) exists, provided the functions $\{f_0, f_1, f_2, \cdots, f_d\}$ satisfy a specific form or representation. If a general solution is computable in polynomial time, then the problem of finding roots can be solved in polynomial time for any ring $\mathbb{R} \subset \mathbb{C}$. The issue is being able to set up the algebraic systems in polynomial time and then solve them in polynomial time, which brings us to the obstacle of determining the formulation of these functions. To inductively find the formula, we must expand the product, and so the time to find each function would be equivalent to the time it takes to expand each product but then after gather the like terms. From the induction proof of Theorem \ref{Thm6}, we know they can be determined inductively; however, one may know the map beforehand and be able to determine if more effectively if it exists regardless of the programmer. Thus, we should be more concerned with the computation of the general solution and how long it will take to find the solution $\textbf{r}$ to $\biggl(f_0(\textbf{r}), f_1(\textbf{r}), \cdots, f_d(\textbf{r})\biggr) = \biggl(a_0, a_1, \cdots, a_d\biggr)$ as each one of these functions will exist regardless if we can find them within an appropriate amount of time. These functions need to be known so that we can attempt various methods to solve them, such as Newton's method, Levenberg-Marquardt Algorithm, Homotopy Continuation Methods, Relaxation Methods (Gauss-Seidel, Jacobi), and Direct Methods (e.g., Gaussian Elimination, or LU Decomposition).  The complexity of solving the system depends on the method chosen and the problem structure. In general, methods like Newton's or Levenberg-Marquardt would take $O(n^3)$ for dense systems, while iterative methods like Gauss-Seidel or Jacobi might perform better for sparse systems. In the context of solving algebraic systems and analyzing time complexity, particularly with respect to the roots of polynomials, there are several factors to consider:\\

\begin{enumerate}
  \item \textbf{Finding the System (Formulation Time Complexity)} \\
  The first step in solving any algebraic system is determining the system itself. According to Theorem \ref{Thm6}, we know that for a given polynomial $p(x)$, the set of functions $$\{ f_0(\textbf{r}), f_1(\textbf{r}), \dots, f_d(\textbf{r}) \}$$ can be computed by expanding the product of factors representing the polynomial. Expanding the product $\prod_{i=1}^k (x - r_i)^{m_i}$ involves determining the coefficients of each term in the expanded product. For each pair of terms, the multiplication and collecting like terms will have polynomial time complexity, with the number of operations depending on the degree $k$ and multiplicities $m_i$ of the roots. Specifically, expanding $(x - r_i)^{m_i}$ would require $O(m_i)$ operations, and summing over all factors would lead to an overall complexity of $O(k \cdot m_{\text{max}})$, where $m_{\text{max}}$ is the largest multiplicity. After expanding, the functions $f_i(r)$ can be collected by summing terms for each degree. This results in a total formulation time complexity of $O(k \cdot m_{\text{max}})$.

  \item \textbf{Solving the System (Solution Time Complexity)} \\
  Once the system is formulated, the next step is solving for the roots $r$, i.e., solving the system of equations:
  \[
  (f_0(r), f_1(r), \dots, f_d(r)) = (a_0, a_1, \dots, a_d)
  \]

  \textbf{Time Complexity of Solving:} The time complexity of solving this system depends on the method employed and the structure of the system. Here are some common methods:

  \begin{itemize}
    \item \textbf{Newton's Method:} Newton's method involves iterating over an initial guess to find roots. For each iteration, computing the Jacobian matrix and solving the linear system requires $O(n^3)$ time for $n$-dimensional systems. The method typically converges quadratically if a good initial guess is chosen \cite{newt1, newt2, newt3}.
    \item \textbf{Levenberg-Marquardt Algorithm:} This method is used for non-linear least squares problems \cite{method1, method2, method3, method4} and generally converges faster than Newton's method for certain types of systems. The time complexity is typically $O(n^3)$, similar to Newton's method, depending on the number of variables and the condition number of the system.
    \item \textbf{Homotopy Continuation Methods:} These methods work by transforming the system into a simpler one whose solution is known, and then continuously deforming the solution to the original system \cite{homotopy1, homotopy2, homotopy3, homotopy4, homotopy5}. The complexity of these methods is generally higher due to the need for many intermediate systems to be solved.
    \item \textbf{Relaxation Methods (e.g., Gauss-Seidel, Jacobi):} These iterative methods are useful for solving systems where the matrix is sparse. Their complexity depends on the number of iterations and the matrix structure but typically have $O(n^2)$ complexity per iteration \cite{cont1, cont2, cont3, cont4, cont5, cont6}.
    \item \textbf{Direct Methods (Gaussian Elimination, LU Decomposition):} These methods \cite{elim1, elim2, elim3, elim4} are deterministic and solve systems exactly \cite{LU1,LU2,LU3, LU4, LU5}, typically requiring $O(n^3)$ time for a general system of size $n$. However, they are less efficient for large sparse systems compared to iterative methods.
  \end{itemize}

  \item \textbf{Combining the Time Complexity} \\
  The total time complexity of the problem can be viewed as a combination of the time to formulate the system and the time to solve it. Based on the above analysis:

  \begin{itemize}
    \item Formulation Time Complexity: $O(k \cdot m_{\text{max}})$
    \item Solution Time Complexity: $O(n^3)$ (for methods like Newton’s or LU Decomposition, assuming dense systems)
  \end{itemize}

  Thus, the total time complexity of solving the algebraic system, considering both finding and solving, would be $O(k \cdot m_{\text{max}} + n^3)$, where $k$ is related to the number of roots, $m_{\text{max}}$ to the multiplicities, and $n$ to the size of the system being solved.

  \item \textbf{P vs NP and Algebraic Systems} \\
  The discussion around \textbf{P} vs $\textbf{Np}$ in this context relates to the general question of whether the problem of finding roots in an arbitrary ring (or even in $\mathbb{C}$) can be solved in polynomial time. If a general solution to the algebraic system (finding roots) could be computed in polynomial time, then it would have profound implications for many problems in algebra and computational complexity. However, the time complexity of formulating the system (which is inherently related to the degree and multiplicities of the roots) and solving it via known methods suggests that the problem is unlikely to be solvable in polynomial time for all cases. This is because solving polynomial systems in the general case is known to be hard, and methods like Newton's method, which work in practice, can still suffer from exponential time complexity in some cases due to the condition number or the need for many iterations.
\end{enumerate}

\section{Summary and Future Work}
In this work, we studied the time complexity of the sliding tile puzzle problem using a matrix representation and the roots of polynomials in $\mathbb{C}[x]$, applying various combinatorial and algebraic methods. We showed that the set of solvable sliding tile configurations is smaller than previously thought, as the solution space is polynomial in relation to the matrix size, suggesting that solving sliding tiles may be easier than expected. The difficulty lies in navigating the solution space and efficiently determining solvable configurations. For polynomial roots, determining the roots efficiently depends on computing the general solution of the algebraic system we developed. Thus, we demonstrate that the time complexity of both problems is tied to the ability to compute the general solution, which differs significantly in each case. Future research will focus on developing more efficient algorithms for determining solvability in higher-dimensional sliding tile puzzles, particularly for larger matrices where computational complexity increases exponentially. We aim to expand the study beyond two- and three-dimensional puzzles to explore more complex, higher-dimensional configurations and their solvability conditions. Additionally, we plan to investigate the polynomial time complexity of sliding tile problems and compare them with other complex computational problems, potentially revealing new relationships with polynomial ring theory. We will also explore heuristic or AI-based approaches to efficiently navigate the solution space and identify solvable configurations without exhaustive search. Finally, we aim to explore how sliding tile properties could be applied in cryptographic algorithms to enhance security and efficiency.
 \newpage

\bibliographystyle{amsplain}

\end{document}